\documentclass{amsart}

\usepackage{palatino}
\usepackage{mathrsfs}
\usepackage{latexsym}

\newcommand{\GL}{\operatorname{GL}}

\newcommand{\PGL}{\operatorname{PGL}}
\newcommand{\SL}{\operatorname{SL}}

\newcommand{\Der}{\operatorname{Der}}
\newcommand{\End}{\operatorname{End}}
\newcommand{\tensor}{\otimes}

\newcommand{\mm}{\mathfrak{m}}

\newcommand\Z{\mathbf{Z}}    %        integers
\newcommand\F{\mathbf{F}}    %        finite fields
\newcommand\Q{\mathbf{Q}}    %        rational field
\newcommand\G{\mathbf{G}}    %        algebraic groups

\newcommand\Aff{\mathbf{A}}      %        affine space
\newcommand\Lie{\operatorname{Lie}}
\newcommand\Ad{\operatorname{Ad}}
\newcommand\ad{\operatorname{ad}}

\newcommand\lie[1]{\mathfrak{#1}}
\newcommand\glie{\lie{g}}

\newcommand\blie{\lie{b}}

\newcommand\plie{\lie{p}}

\newcommand{\hlie}{\lie{h}}
\newcommand{\ulie}{\lie{u}}
\newcommand{\vlie}{\lie{v}}

\newcommand{\clie}{\lie{c}}

\newcommand{\zlie}{\lie{z}}

\newcommand{\normal}{\lhd}

\newcommand{\congruent}{\equiv} 
\newcommand{\iso}{\simeq}

\newtheorem{theorem}{Theorem}
\newtheorem{prop}[theorem]{Proposition}

\newtheorem{lem}[theorem]{Lemma}

\theoremstyle{remark}
\newtheorem{rem}[theorem]{Remark}

\newcommand{\Zp}{\Z_{(p)}}
\newcommand{\e}{\varepsilon}

\newcommand{\A}{\mathscr{A}}

\newcommand{\UU}{\mathcal{U}}
\newcommand{\VV}{\mathcal{V}}
\newcommand{\NN}{\mathcal{N}}
\newcommand{\MM}{\mathcal{M}}

\newcommand{\Int}{{\operatorname{Int}}}
\newcommand{\rank}{{\operatorname{rank}}}

\DeclareMathAlphabet{\matheucal}{U}{eus}{m}{n}
\newcommand{\Orbit}{\matheucal{O}}
\newcommand{\Class}{\matheucal{C}}

\newcommand{\algc}[1]{\overline{#1}}

\numberwithin{equation}{subsection}

\title{Sub-principal homomorphisms in positive characteristic}
\author{George J. McNinch}
\thanks{Supported in part by a grant from the National Science Foundation}
\address{Department of Mathematics \\
         Room 255 Hurley Building \\
         University of Notre Dame \\
         Notre Dame, Indiana 46556-5683 \\
         USA}
\email{mcninch.1@nd.edu}
\date{June 10, 2002}

\begin{document}
\bibliographystyle{amsalpha}

\begin{abstract}
  Let $G$ be a reductive group over an algebraically closed field of
  characteristic $p$, and let $u \in G$ be a unipotent element of
  order $p$. Suppose that $p$ is a good prime for $G$. We show in this
  paper that there is a homomorphism $\phi:\SL_{2/k} \to G$ whose
  image contains $u$. This result was first obtained by D. Testerman
  (J. Algebra, 1995) using case considerations for each type of simple
  group (and using, in some cases, computer calculations with explicit
  representatives for the unipotent orbits).
  
  The proof we give is free of case considerations (except in its
  dependence on the Bala-Carter theorem). Our construction of $\phi$
  generalizes the construction of a principal homomorphism made by
  J.-P.  Serre in (Invent. Math. 1996); in particular, $\phi$ is
  obtained by reduction modulo $\mathfrak{p}$ from a homomorphism of
  group schemes over a valuation ring $\A$ in a number field.  This
  permits us to show moreover that the weight spaces of a maximal
  torus of $\phi(\SL_{2/k})$ on $\Lie(G)$ are ``the same as in
  characteristic 0''; the existence of a $\phi$ with this property was
  previously obtained, again using case considerations, by Lawther and
  Testerman (Memoirs AMS, 1999) and has been applied in some recent
  work of G.  Seitz (Invent. Math. 2000).
\end{abstract}

\maketitle

%\tableofcontents

\section{Introduction}
\label{sec:intro}

Let $G = G_{/k}$ be a connected reductive algebraic group over an algebraically
closed field $k$ of characteristic $p > 0$.  It is the main goal of this
note to give another proof of the following theorem:

\begin{theorem} (Testerman \cite{testerman}) 
  \label{theorem:testerman}
  Suppose that $p$ is a good prime for $G$.  If $u \in G$ is unipotent
  and has order $p$, then there is a homomorphism $\phi:\SL_{2/k} \to G$
  with $u$ in its image.
\end{theorem}

One might regard Theorem \ref{theorem:testerman} as a group analogue
of the Jacobson-Morozov theorem for Lie algebras. If one considers
instead any field $E$ of characteristic 0, a reductive group $G_{/E}$
over $E$, and $u \in G_{/E}$ an $E$-rational unipotent element, one
may write $u = \exp(X)$ for a nilpotent $X \in \Lie(G_{/E})$; from the
Jacobson-Morozov theorem for $\Lie(G_{/E})$ one deduces a homomorphism
$\SL_{2/E} \to G_{/E}$ over $E$ with $u$ in its image.

Testerman's original proof of Theorem \ref{theorem:testerman} used case
considerations for each type of simple group (and used, in some cases,
computer calculations with explicit unipotent class representatives
known from the work of Mizuno).  Our proof for the most part avoids
case considerations (except that it depends on Pommerening's proof of
the Bala-Carter theorem in good characteristic).  We will exploit a
weak Jacobson-Morozov-type result for an integral form of the Lie
algebra of $G$. We obtain $\phi$ as a suitable $G$-conjugate of the
reduction mod $\mathfrak{p}$ of a homomorphism of group schemes
$\phi_{/\A}:\SL_{2/\A} \to G_{/\A}$, where $\A$ is a valuation ring in
a number field.

When $u$ is regular unipotent and has order $p$, the theorem yields a
\emph{principal homomorphism}; see \cite[\S 2]{serre:PSLp}. The
argument we give specializes in the regular case to the proof of
\emph{loc.  cit.}  Proposition 2.  

The fact that we obtain a homomorphism of group schemes over
$\A$ permits us to prove a more precise version of Theorem
\ref{theorem:testerman}, which was first obtained by Lawther and
Testerman. 

In good characteristic, one can associate to a nilpotent element $X
\in \glie$ a cocharacter $\nu$, which is well defined up to
conjugation by $C_G^o(X)$ (the connected centralizer); see Proposition
\ref{prop:associated-cocharacter}.  Moreover, we may find a
$G$-equivariant homeomorphism $\e$ from the nilpotent variety to the
unipotent variety; see Proposition \ref{prop:springer-homeo}. We say
that $\e$ is a Springer homeomorphism.

We say that a homomorphism $\phi:\SL_{2/k} \to G$ is
\emph{sub-principal}
  \footnote{  To explain this terminology, note that the main result of
  \cite{Lawther-Testerman} shows (under some conditions on $p$ which
  are slightly more restrictive than ``$p$ is good'') that there is a
  unique conjugacy class of principal homomorphisms $\phi$ with $d
  \phi \not = 0$, and that for each unipotent class there is a unique
  conjugacy class of sub-principal homomorphisms (in the sense defined
  in this paper). Thus the classes of sub-principal homomorphisms are
  in some sense analogous to the class of principal homomorphisms
  (they are precisely the classes which ``come from characteristic
  0'').}
if the restriction of $\phi$
to a maximal torus of $\SL_{2/k}$ is a cocharacter associated to
some non-0 nilpotent element $X$ in the image of $d\phi$, and if
 $\e(X)$ is conjugate to a unipotent element in the image of $\phi$.
 
 Note that by Lemma \ref{lem:iso-respects-parabolics} and the
 Bala-Carter Theorem (Proposition \ref{prop:bala-carter}), if $\e'$ is
 another Springer homeomorphism, then $\e(X)$ is conjugate to
 $\e'(X)$. Thus the notion of a sub-principal homomorphism is
 independent of this choice.

We may now state the more precise form of Theorem
\ref{theorem:testerman}:
\begin{theorem}
\label{theorem:LawtherTesterman}
(Lawther and Testerman \cite[Theorem 4.2]{Lawther-Testerman}) With the
assumptions of Theorem \ref{theorem:testerman}, there is a
homomorphism $\phi:\SL_{2/k} \to G$ such that $u$ is in the image of
$\phi$,  the restriction of $\phi$ to a maximal
torus of $\SL_{2/k}$ is a co-character associated to 
some nilpotent $0 \not = X \in \text{Image}(d\phi)$, and 
$\e(X)$ is conjugate to  $u$.
Thus $\phi$ is a \emph{sub-principal} homomorphism.
\end{theorem}

In the language used by Lawther and Testerman
\cite{Lawther-Testerman}, the theorem yields an $A_1$ subgroup of $G$
whose labeled Dynkin diagram is the same as the labeled diagram
of $u$. Indeed, the labeled diagram of the $A_1$ subgroup is obtained
by choosing a maximal torus $T_0$ of $\SL_{2/k}$ and maximal torus $T$
of $G$ containing $\phi(T_0)$. The homomorphism $\mu=\phi_{\mid
  T_0}:T_0 \iso \G_m \to G$ is then a cocharacter.  For a suitable
choice of Borel subgroup $B$ containing $T$ (equivalently: a suitable
choice of positive roots) the values $\langle \alpha,\mu \rangle$ at
the simple roots in $X^*(T)$ are non-negative and constitute the
labels on the Dynkin diagram. One checks that these labels are
independent of the choices made; see \cite[\S 7.6]{hum-conjugacy}.
Similarly, the labels on the diagram of $u$ are the non-negative
integers $\langle \alpha,\nu \rangle$ where $\nu$ is a co-character
associated to $u$ (where again $T$ and $B$ are suitably chosen).

Now let $E$ be an algebraically closed field of characteristic 0, and
let $G_{/E}$ be a reductive group over $E$ with the same root datum as
$G$.  There is a bijection between unipotent classes
in $G$ and unipotent classes in $G_{/E}$ which preserves labeled diagrams. The Dynkin-Kostant classification of
nilpotent orbits in characteristic 0 implies: If $\phi_{/E}:\SL_{2/E}
\to G_{/E}$ is any homomorphism with the unipotent element $v$ in its image, the
labeled diagram of $\phi_{/E}$ coincides with that of $v$.  Thus
Theorem \ref{theorem:LawtherTesterman} yields a
homomorphism $\phi:\SL_{2/k} \to G$ whose image contains $u$ and for
which the weights of a maximal torus of $\SL_{2/k}$ on $\Lie(G)$ are
``the same as in characteristic 0.'' (We refer the reader to the
extensive tables in \cite{Lawther-Testerman} to see that for some $u$
there are homomorphisms $\phi:\SL_{2/k} \to G$ whose image contains $u$,
but whose labeled diagram differs from that of $u$).

We mention that Theorem \ref{theorem:LawtherTesterman} was used by
Seitz in \cite{seitz-unipotent}. In \emph{loc. cit.}, Seitz introduced
the notion of a ``good $A_1$''. In the language above, he calls a
homomorphism $\phi:\SL_2(k) \to G$ \emph{good} if each weight $\lambda
\in \Z$ of the representation $(\Ad \circ \phi,\glie)$ for a maximal
torus of $\SL_2(k)$ satisfies $\lambda \le 2p-2$, where $\glie =
\Lie(G)$.  Seitz proves that for each unipotent $u$ of order $p$,
there is a good homomorphism $\phi$ with $u$ in its image; his proof
of the existence of such a $\phi$ depends in a crucial way on Theorem
\ref{theorem:LawtherTesterman} (when $u$ is distinguished, the
existence of a good $\phi$ is an immediate consequence of Theorem
\ref{theorem:LawtherTesterman} combined with the ``order formula'' of
Testerman which one will find in \cite{testerman} or
\cite[Theorem 1.1]{math.RT/0007056}).

We also obtain a refinement of Theorem \ref{theorem:LawtherTesterman}
for finite fields: if $G$ is defined over a finite field $\F_q$ of
good characteristic and $u$ is $\F_q$-rational, we show in \S
\ref{sub:finite-fields} that $\phi$ may be defined over $\F_q$ as
well.  In the course of our proof, we establish the following result
which may be of independent interest. Let $G$ be defined over an
arbitrary field $F$ (of good characteristic) and let $u$ be an
$F$-rational unipotent element. Suppose that either the orbit of $u$
is separable, or that $F$ is perfect. Then the canonical parabolic
subgroup attached to $u$ is defined over $F$. (The same statement
holds for $F$-rational nilpotent elements).

Finally, we present two appendices. In the first, we give a proof that
in good characteristic, there is always a $G$-equivariant
\emph{homeomorphism} between the nilpotent variety $\NN$ and the
unipotent variety $\UU$ of a reductive group. Of course, in ``very
good'' characteristic, there is an isomorphism of varieties due to
Springer; our argument handles also groups such as
$\operatorname{PGL}_{p/k}$ in characteristic $p$. This simplifies some
of the steps in our proof of Theorem \ref{theorem:LawtherTesterman}.
In the second, we show that the $G$-equivariant isomorphism $\UU \iso
\NN$ defined by Bardsley and Richardson \cite{Bard-Rich-LunaSlice}
respects the $p$-th power operations.

The author would like to thank the Mathematics Institute at the
University of {\AA}rhus for its hospitality during the academic year
2000/2001. Thanks also to Jens C. Jantzen and Donna M. Testerman for
some suggestions on this manuscript. Finally, thanks to the referee
for some useful remarks concerning \S \ref{sub:finite-fields}.

\section{Generalities on reductive groups}

Let $G_{/\Z}$ be a split reductive group scheme over $\Z$. If
$G_{/\Z}$ is moreover semisimple, one may regard $G_{/\Z}$ as a
``Chevalley group scheme'' as in
\cite{borel70:_proper_linear_repres_cheval_group}. Let $\glie_\Z$ be
the Lie algebra.  For any commutative ring $\Lambda$, we put
$G_{/\Lambda} = G_{/\Z} \times_\Z \Lambda$, and $\glie_\Lambda =
\glie_\Z \tensor_\Z \Lambda$.

Let $(X,Y,R,R^\vee)$ be the root datum of $G_{/\Z}$ with respect to a
fixed maximal torus $T_{/\Z}$.  Fix a $\Z$-basis
$\gamma_1,\dots,\gamma_r$ for $Y$, the co-character group of
$T_{/\Z}$. Let $H_{\gamma_i} = d\gamma_i(1) \in \hlie_\Z$, the Lie
algebra of $T_{/\Z}$. The algebra $\glie_\Z$ has a Chevalley basis
$$\{E_\alpha \mid \alpha \in R\} \cup
\{H_{\gamma_1},\dots,H_{\gamma_r}\}.$$ We have $\hlie_\Z = \bigoplus_i
\Z H_{\gamma_i}$, and $\blie_\Z =
\hlie_\Z \oplus \bigoplus_{\alpha \in R^+} \Z E_\alpha$ is a Borel
subalgebra of $\glie_\Z$.

\subsection{Good primes}
\label{sub:good-p}

Recall the notion of a good prime for the root system $R$ (or for the
root datum $(X,Y,R,R^\vee)$, it is the same).

For the indecomposable root systems, a prime is bad (= not good) only
in the following situations: $2$ is bad unless $R$ is of type $A$, $3$
is bad if $R$ is of type $E,F$ or $G$, and $5$ is bad if $R=E_8$. 

For general $R$, $p$ is good for $R$ provided that it is good
for each indecomposable component of $R$.

\subsection{Parabolic subalgebras}
\label{sub:parabolics}

If $S$ denotes the simple roots in $R^+$, any subset $I \subset S$
determines a subroot system $R_I$ in a well-known way, and
hence a \emph{standard parabolic subalgebra}
\begin{equation*}
  \plie(I)_\Z = \blie_\Z \oplus 
  \bigoplus_{\alpha \in R_I^+} \Z E_{-\alpha}.
\end{equation*}

Consider the function $f:\Z R \to \Z$ which satisfies
\begin{equation}
  \label{eq:almost-cocharacter}
  f(\alpha) = \left \{
 \begin{matrix}
   2 & \alpha \in S \setminus I \\
   0 & \alpha\in I
 \end{matrix} \right.
\end{equation}
We may regard $f$ as a co-character of the adjoint group, so that
$\glie_\Z$ becomes a module via $f$ for $\G_{m/\Z}$; as such, it is
the direct sum of its weight spaces. Thus, we have $\glie_\Z =
\bigoplus_{i \in \Z} \glie_\Z(i)$, where
\begin{equation*}
  \glie_\Z(0) = \hlie_\Z \oplus \bigoplus_{f(\alpha) = 0} 
  \Z E_\alpha, \quad  \text{and} \quad
  \glie_\Z(i) = \bigoplus_{f(\alpha) = i} 
  \Z E_\alpha \quad \text{for $i \not = 0$}.
\end{equation*}
We obtain the original parabolic algebra as $\plie(I)_\Z =
\bigoplus_{i \ge 0} \glie_\Z(i)$.  The opposite parabolic subalgebra
is $\plie(I)^-_\Z = \bigoplus_{i \le 0} \glie_\Z(i)$.  We put
$\ulie(I)_\Z = \bigoplus_{i>0} \glie_\Z(i)$ and $\ulie(I)^-_\Z =
\bigoplus_{i<0} \glie_\Z(i)$.

There are ``group scheme versions'' of each of these constructions:
i.e. there are parabolic subgroup schemes $P(I)_{/\Z}$ and
$P(I)^-_{/\Z}$ with respective subgroup schemes $U(I)_{/\Z}$ and
$U(I)^-_{/\Z}$.

\subsection{Distinguished nilpotents and parabolics}

Let $k$ be an algebraically closed field with characteristic $p \ge
0$; in this section and the next we write $G=G_{/k}$ and $\glie =
\glie_k$. We suppose that $p$ is good for $G$.  A nilpotent element $X
\in \glie$, respectively a unipotent element $u \in G$, is said to be
\emph{distinguished} if the connected center of $G$ is a maximal torus
of $C_G(X)$, respectively $C_G(u)$. A parabolic subalgebra $\plie \subset \glie$ is
called \emph{distinguished} if
\begin{equation*}
  \dim \plie / \ulie = \dim \ulie /[\ulie,\ulie] + \dim \Lie(Z),
\end{equation*}
where $\ulie$ is the nilradical of $\plie$ and $Z$ is the center of $G$.

Let $\plie_\Z = \plie(I)_\Z$ be a standard parabolic subalgebra of
$\glie_\Z$ as in \S \ref{sub:parabolics}, and let $\ulie_\Z =
\ulie(I)_\Z$. Then $\ulie = \ulie_\Z \tensor_\Z k$ is the nilradical
of $\plie = \plie_Z \tensor_\Z k$.  Let $\glie(i) = \glie_\Z(i)
\tensor_\Z k$ for $i \in \Z$. Since $p$ is good, $\dim_k \glie(2) =
\dim \ulie/[\ulie,\ulie]$; see \cite[Prop.  4.3]{math.RT/0007056} or
the proof of \cite[Prop. 5.8.1]{Carter1}.

Thus the condition that $\plie_k$ be
distinguished is independent of $k$ so long as the characteristic of
$k$ is good; we say that $\plie_\Z$ is a distinguished standard
parabolic subalgebra if $\plie_k$ is distinguished.

When $\plie_\Z$ is distinguished, it
follows from \cite[Lemma 5.2]{jantzen:Nilpotent} that
the map $f:\Z R \to \Z$ of
\eqref{eq:almost-cocharacter} extends uniquely to a linear function
$X^*(T_{/\Z}) \to \Z$ and hence determines a cocharacter $\tau$ of
$T_{/\Z}$ satisfying:
\begin{equation}
  \label{eq:cocharacter}
 \langle \alpha,\tau \rangle = f(\alpha).
\end{equation}
Note that the argument in \emph{loc. cit.} applies to semisimple $G$,
which we may reduce to by considering the derived group of $G$.

\subsection{Richardson orbits and the Bala-Carter Theorem}
\label{section:bala-carter-pom}

Let $k$, $G$, $\glie$ as in the previous section; especially, recall
that $p$ is good. Suppose now that $\plie \subset \glie$ is any
parabolic subalgebra, with nilradical $\ulie$.  There is a unique
parabolic subgroup $P \le G$ with $\Lie(P) = \plie$.  Moreover,
$\ulie$ is the Lie algebra of the unipotent radical $U$ of $P$.  A
theorem of R. Richardson \cite[Theorem 5.3]{hum-conjugacy} says that
there is a nilpotent $G$-orbit $\Orbit \subset \glie$ with the
property that $\Orbit \cap \ulie$ is an open $P$-orbit. Similarly,
there is a unipotent class $\Class$ in $G$ with the property that
$\Class \cap U$ is an open $P$-orbit.  By a Richardson element, we
mean an orbit representative for $\Class$ or $\Orbit$ lying in $U$
respectively $\ulie$. The orbits $\Orbit$ and $\Class$ are known as
the Richardson orbits associated with $\plie$ (or with $P$).

By a Levi subgroup of $G$, we mean a Levi factor of a parabolic
subgroup.

\begin{prop} (Bala-Carter, 
              Pommerening \cite{PommereningI,PommereningII})
             \label{prop:bala-carter}
  \begin{enumerate}
  \item Consider the collection of all pairs $(L,\Orbit)$ consisting
    of a Levi subgroup of $G$ and a distinguished nilpotent orbit in
    $\Lie(L)$. Then the map which associates to $(L,\Orbit)$ the
    $G$-orbit $\Ad(G)\Orbit$ defines a bijection between the set of
    $G$-orbits of pairs $(L,\Orbit)$ and the nilpotent $G$-orbits in
    $\glie$.
  \item Associate to each distinguished parabolic subalgebra its
    nilpotent Richardson orbit. Then this map defines a bijection
    between the conjugacy classes of distinguished parabolic
    subalgebras and the distinguished nilpotent orbits in $\glie$.
  \end{enumerate}
\end{prop}
Note that (1) holds with no assumption on $p$, but that (2) requires
$p$ to be good.

A cocharacter $\nu:\G_m \to G$
is said to be associated to a nilpotent element $X \in \glie$
provided that $\Ad(\nu(t))X = t^2X$ for all $t \in \G_m$,
and that $\nu$ takes values in the derived group of some Levi
subgroup $L$ for which $X \in \Lie(L)$ is distinguished.

\begin{rem}
  \label{rem:levi-assoc-cocharacter}
  Let $L$ be a Levi subgroup of a parabolic in $G$, and let $X$ be a
  nilpotent element in $\Lie(L)$. If $\iota:L \to G$ is the inclusion
  map, then a cocharacter $\tau$ of $L$ is associated to $X$ (with
  respect to $L$) if and only if $\iota \circ \tau$ is associated to
  $X$ (with respect to $G$).  If $\e:\NN \to \UU$ is a $G$-equivariant
  homeomorphism, then Lemma \ref{lem:iso-respects-parabolics} (in the
  appendix) shows that $\e$ restricts to a suitable $L$-equivariant
  homeomorphism. It follows that $\phi:\SL_{2/k} \to L$ is a
  sub-principal homomorphism (for $L$) if and only if $\iota \circ
  \phi$ is a sub-principal homomorphism (for $G$).
\end{rem}

\begin{rem}
  \label{rem:cochar-and-isogeny} Let $\pi:G \to G'$ be a central
  isogeny of reductive groups.  According to Lemma
  \ref{lem:central-isogeny-homeo}, $d\pi$ restricts to a
  $G$-equivariant homeomorphism between the respective nilpotent
  varieties. If $L$ is a Levi subgroup of $G$ then $\pi(L) = L'$ is a
  Levi subgroup of $G'$, and it follows  that $X$ is distinguished in
  $\Lie(L)$ if and only if $d\pi(X)$ is distinguished in
  $\Lie(L')$. So if $X$ is nilpotent in $\Lie(G)$, and if
  $\phi$ is a co-character associated to $X$, then $\pi \circ \phi$ is
  a co-character associated to $d\pi(X)$.
\end{rem}

We record the following:
\begin{prop}
  \label{prop:associated-cocharacter}
  Suppose the characteristic is good for $G$.  There is a cocharacter
  associated to any nilpotent element $X$.  Moreover, any two
  cocharacters associated to $X$ are conjugate by an element of
  $C^o_G(X)$.
\end{prop}

\begin{proof}
  \cite[Lemma 5.3]{jantzen:Nilpotent}
\end{proof}

\begin{prop}
  \label{prop:canonical-parabolic}
  Suppose the characteristic is good for $G$.  Let $X$ be nilpotent,
  and let $\nu$ be a cocharacter associated to $X$. Consider the
  parabolic subalgebra $\plie = \bigoplus_{i \ge 0} \glie(i)$, and let
  $P$ be the corresponding parabolic subgroup of $G$.
  \begin{enumerate}
  \item $C_G(X) < P$.
  \item $P$ and $\plie$ depend only on $X$.
  \end{enumerate}
\end{prop}

\begin{proof}
  \cite[Prop. 5.9]{jantzen:Nilpotent}
\end{proof}

The subalgebra $\plie$ is known as the \emph{canonical} (or
\emph{Jacobson-Morozov}) parabolic determined by $X$.

\subsection{$\lie{sl}_2$ triples}

Let $\Lambda$ be an integral domain. The data $(0,0,0) \not = (X,Y,H)
\in \glie_\Lambda \times \glie_\Lambda \times \glie_\Lambda$ is called
an $\lie{sl}_2$ triple (over $\Lambda$) if the $\Lambda$-linear map
$d\phi:\lie{sl}_2(\Lambda) \to \glie_\Lambda$ given by $
\begin{pmatrix}
  0 & 1 \\
  0 & 0
\end{pmatrix} \mapsto X$,
$
\begin{pmatrix}
  0 & 0 \\
  1 & 0 
\end{pmatrix} \mapsto Y$,
and $
\begin{pmatrix}
  1 & 0 \\
  0 & -1
\end{pmatrix} \mapsto H$ is an injective Lie algebra homomorphism.
To check that $d\phi$ is a homomorphism, one only needs to see that
$[X,Y] = H$, $[H,X] = 2X$ and $[H,Y]=-2Y$. If the characteristic
of $\Lambda$ is not 2, injectivity is immediate.

\begin{prop}
\label{prop:sl2-triple-char0}
  Let $F$ be a field of characteristic 0, and let $(X,Y,H)$
be an $\lie{sl}_2$ triple in $\glie_F$. Then
there is a unique homomorphism $\phi:\SL_{2/F} \to G_{/F}$
whose tangent map is the Lie algebra homomorphism $d\phi$ associated
with the triple $(X,Y,H)$.
\end{prop}

\begin{proof}
  By \cite[Prop. 7.1]{math.RT/0007056}, there is a homomorphism $t
  \mapsto \exp(tX) :\G_{a/F} \to G_{/F}$ for which $\rho(\exp(tX)) =
  \exp(d\rho(tX))$ for every rational representation $\rho$ of
  $G_{/F}$.
  
  Let $(\rho,V)$ be a faithful $F$-representation of $G_{/F}$. Then
  $d\rho$ restricts to a representation of the $\lie{sl}_2$ triple.
  The Chevalley group construction \cite{Steinberg} applied to the Lie
  algebra $\lie{sl}_{2/F}$ and the representation $(d\rho,V)$ gives a
  homomorphism $\phi:\SL_{2/F} \to \GL(V)$ which maps the upper
  triangular subgroup to the image of $t \mapsto \exp(d\rho(tX))$ and
  the lower triangular subgroup to the image of $s \mapsto \exp(
  d\rho(sY))$.  Now \cite[Prop.  6.12]{Bor1} implies that $\phi$ takes
  values in $G_{/F}$ and is unique. It is clear by construction that
  $d\phi$ has the desired form.
\end{proof}

\section{The main result}

In this section, we suppose that $k$ is an algebraic closure of the
finite field $\F_p$.  The split reductive group scheme $G_{/\Z}$ is as
in the previous section; we now suppose that $G_{/\Z}$ is (split)
\emph{semisimple} and \emph{simply connected}. In particular, $G=G_{/k}$
is simply connected. Note that we reserve the undecorated notations $G$,
$\glie$, etc. for the objects over $k$. We suppose $p$ to be good for
$G$.

Fix an algebraic closure $E=\algc{\Q}$ of the rational field; in what
follows, we regard all finite extensions of $\Q$ as subextensions of
$E/\Q$.  For $F$ a finite extension of $\Q$ and $\A \subset F$ a
valuation ring whose residue field has characteristic $p$, an element
$X \in \glie_\A$ determines elements $X_F = X \tensor 1_F \in \glie_F$
$X_E \in \glie_E$, and $X_k \in \glie = \glie_k$ (note that $X_k$
actually depends on  the embedding of the residue field of
$\A$ in $k$; the particular choice is not important).

\subsection{$\lie{sl}_2$-triples over integers}

%\begin{lem}
%  \label{lem:inf-centralizer-lemma}
%  If $X \in \glie$ is a distinguished nilpotent element satisfying
%  $X^{[p]} = 0$, then $\clie_\glie(X) \subset \plie$, where $\plie$ is
%  the canonical parabolic subalgebra attached to $X$.
%\end{lem}

%\begin{proof}
%  Let $\pi:\tilde G \to G$ be a central isogeny with $\tilde G$ a
%  direct product of simply connected, quasisimple groups $G_1 \times
%  \cdots \times G_r$. \cite[Prop.  2.6]{jantzen:Nilpotent} (or Lemma
%  \ref{lem:central-isogeny-homeo} in the appendix) shows that there is
%  a unique nilpotent element $\tilde X \in d\pi^{-1}(X)$.  It is clear
%  that $\tilde X$ is distinguished in $\tilde \glie = \Lie(\tilde G)$.
  
%  For $i=1,\dots,r$, let $p_i:\tilde G \to G_i$ be the projection.
%  Then $X_i = p_i(\tilde X)$ is distinguished in $\glie_i = \Lie(G_i)$
%  for each $i$. If the root system of $G_i$ is of type $A_{n(i)}$ for
%  $1 \le i \le r$, then $G_i \iso \SL_{n(i)}(k)$. Thus $X_i$ is
%  regular nilpotent. The condition $X_i^{[p]} = 0$ then implies that
%  $n(i) \ge p$. The survey in \cite[0.13]{hum-conjugacy} now shows
%  (since the characteristic $p$ is good for $G$) that $(\Ad,\glie_i)$
%  is a simple $G_i$-module for each $1 \le i \le r$; it follows that
%  $d\pi$ is bijective. Thus, we may suppose that $G=\tilde G$. The
%  assertion in that case follows from the main result of
%  \cite{spalt-transverse}.
%\end{proof}

We require the following result due to Spaltenstein.  The simple
connectivity hypothesis is unnecessary when $G$ has no simple factors
of type $A_n$.  Some hypothesis is necessary, though, since the
conclusion of the following lemma is not valid for example when $p=2$
and $G$ is the adjoint group $\PGL_2(k)$.

\begin{lem}
  \label{lem:inf-centralizer-lemma}
  If $X \in \glie$ is a distinguished nilpotent element, then
  $\clie_\glie(X) \subset \plie$, where $\plie$ is the canonical
  parabolic subalgebra attached to $X$.
\end{lem}

\begin{proof}
  Since $G$ is simply connected, it is a direct product of simply
  connected, quasi-simple groups $G \iso G_1 \times \cdots \times G_r$.
  For $i=1,\dots,r$, let $p_i:G \to G_i$ be the projection.  Then $X_i
  = dp_i(X)$ is distinguished in $\glie_i = \Lie(G_i)$ for each $i$.
  Since $\clie_\glie(X) = \bigoplus_i \clie_{\glie_i}(X_i)$, and since
  $\plie = \bigoplus_i \plie \cap \glie_i$, it suffices to assume that
  $G$ is quasi-simple.  
  
  If the root system of $G$ is not of type $A_n$, the assertion
  follows from the main result of Spaltenstein in
  \cite{spalt-transverse}.  Otherwise, $G \iso \SL_n(k)$, and $X$ is a
  regular nilpotent element in $\glie = \lie{sl}_n(k)$. In this case,
  $\plie = \blie$ is the Lie algebra of a Borel subgroup, and the
  claim follows from a direct computation.
\end{proof}

\begin{lem}
  \label{lem:integral-sl2}
  Let $\plie_\Z \subset \glie_\Z$ be a distinguished standard
  parabolic, and $\tau \in X_*(T_{/\Z})$ the corresponding co-character
  as in \eqref{eq:cocharacter}.  Suppose that a Richardson element $X
  \in \ulie_k$ satisfies $X^{[p]} = 0$.

  \begin{enumerate}
  \item There is a finite field extension $F \supset \Q$, a valuation
    ring $\A \subset F$ (whose residue field we embed in $k$), and an
    element $X \in \glie_\A(2)$ such that 
    \begin{enumerate}
    \item[(a)] $X_k \in \glie(2)$ is
    a Richardson element for $\plie$,
    \item[(b)] $X_E \in \glie_E(2)$
    is a Richardson element for $\plie_E$.
    \end{enumerate}
  \item Put $H = d\tau(1) \in \hlie_\Z \subset \glie_\Z(0)$. There is
    a unique element $Y \in \glie_\A(-2)$ such that $(X,Y,H)$ is an
    $\lie{sl}_2$-triple over $\A$.
  \end{enumerate}
\end{lem}

\begin{proof}
  The assertion (1) is elementary; a proof is written down in
  \cite[Lemma 5.2]{math.RT/0007056}.  Since $H = d\tau(1)$, (2) will
  follow provided that we find $Y \in \glie_\A(-2)$ with $[X,Y] = H$.
  Since $\plie$ is distinguished and  $G_{/\Z}$ is semisimple,
  we have $$\rank_\Z \ \glie_\Z(-2) = \rank_\Z \ \glie_\Z(0).$$  Lemma
  \ref{lem:inf-centralizer-lemma} shows that the centralizer in $\glie$
  of $X_k$ is contained in $\plie$. This implies that
  $\ad(X_k):\glie(-2) \to \glie(0)$ is injective, and is therefore
  a linear isomorphism.  Thus $\ad(X):\glie_\A(-2) \to \glie_\A(0)$ is
  also bijective and (2) follows.
\end{proof}

\subsection{Exponential isomorphism}
\label{sub:exponential-iso}

\begin{lem}
\label{lem:exponential-iso}
  Let $\plie_\Z$ be a distinguished standard parabolic subalgebra of
  $\glie_\Z$, and let $\ulie_\Z$ be as
  in \ref{sub:parabolics}. Suppose that a Richardson element $X \in
  \ulie_k$ satisfies $X^{[p]} = 0$. Then the exponential isomorphism
  $\exp:\ulie_\Q \to U_\Q$ is defined over $\Zp$.
\end{lem}

\begin{proof}
  Recall that $\ulie_\Z$ is the Lie algebra of the group scheme
  $U_\Z$; $U_\Q$ is obtained by base change.  Since $X^{[p]} = 0$,
  \cite[Theorem 5.4]{math.RT/0007056} shows that that $\glie_\Z(i) =
  0$ for all $i \ge 2p$ (for the grading induced by the cocharacter
  $\tau$ in \eqref{eq:cocharacter}). It now follows 
    %%%%%%%%
  \footnote{
    In \cite[\S 4.4]{math.RT/0007056}, $n(P)$ is defined as the least
    $n \ge 0$ with $\glie(2n) = 0$ for the grading induced by $\tau$
    as in \eqref{eq:cocharacter}.  Write $c(\ulie)$ for the nilpotence
    class of $\ulie$, and $c(U)$ for that of $U$.  \cite[Prop.
    4.4]{math.RT/0007056} erroneously asserts that $c(\ulie)$, $c(U)$
    and $n(P)$ coincide; in fact the given proof shows that $n(P)-1 =
    c(\ulie) = c(U)$.  In the above situation, we therefore see that
    $p \ge n(P) > c(\ulie)$ whence the claim.  This error in
    \cite[Prop.  4.4]{math.RT/0007056} led to a flawed statement of
    \cite[Theorem 1.1]{math.RT/0007056}; a correct statement is
    obtained by taking $n(P)$ to be $c(\ulie)+1$ (rather than
    $c(\ulie))$.
    }
  %%%%%%%%
  that the nilpotence class of the Lie algebra $\ulie_k$ is $< p$.
  The lemma now follows from \cite[Prop. 5.1]{seitz-unipotent}.
\end{proof}

\begin{lem}
  \label{lem:rational-exponential-iso}  
  Let $F$ be an arbitrary subfield of $k$, and suppose that $G$ is
  defined (not necessarily split) over $F$.  Let $\plie$ be a
  parabolic subalgebra defined over $F$, and suppose that a Richardson
  element $X$ of its nilradical $\ulie$ satisfies $X^{[p]} = 0$. Then
  there is a unique $P$-equivariant isomorphism of varieties
  $\exp:\ulie \to U$  whose tangent map is
  the identity. Moreover, $\exp$ is defined over $F$.
\end{lem}

\begin{proof}
  As in the proof of the previous lemma, the hypothesis guarantees
  that the nilpotence class of $\ulie$ is $< p$; the result then
  follows from \cite[Proposition 5.2]{seitz-unipotent}.
\end{proof}

\subsection{The main theorem}

It is:
\begin{theorem}
  \label{theorem:main}
  Let $\plie_\Z \subset \glie_\Z$ be a distinguished standard
  parabolic, and $\tau \in X_*(T_{/\Z})$ the corresponding
  co-character as in \eqref{eq:cocharacter}. Suppose that a Richardson
  element $Y \in \ulie_k$ satisfies $Y^{[p]} = 0$.

  Choose the number field $F$, valuation ring $\A$, and element $X \in
  \glie_\A(2)$ as in Lemma \ref{lem:integral-sl2}(1), and let
  $(X,Y,H)$ be an $\lie{sl}_2$-triple over $\A$ as in that
  lemma.
  \begin{enumerate}
  \item Let $\phi_{/F}:\SL_{2/F} \to G_{/F}$ be the homomorphism
    determined by the $\lie{sl}_2$ triple $(X_F,Y_F,H_F)$ as in
    Proposition \ref{prop:sl2-triple-char0}.  Then $\phi_{/F}$ is
    defined over $\A$; i.e. there is a homomorphism of group schemes
    $\phi_{/\A}:\SL_{2/\A} \to G_{/\A}$ from which $\phi_{/F}$ arises
    by scalar extension.
  \item If $\phi_{/k}:\SL_{2/k} \to G$ denotes the map obtained by base
    change from $\phi_{/\A}$, then:
    \begin{enumerate}
    \item The image of $\phi_{/k}$ meets the open $P$-orbit on $U$,
      where $U$ is the unipotent radical of the parabolic subgroup $P$
      corresponding to $\plie$.
    \item The Richardson element $X_k$ is in the image of
      $d\phi_{/k}$.
    \item The restriction of $\phi_{/k}$ to a suitable maximal torus
      of $\SL_{2/k}$ coincides with the co-character $\tau$; it is a
      co-character associated to $X_k$.
    \end{enumerate}
  \end{enumerate}
\end{theorem}

\begin{proof}
  Lemma \ref{lem:exponential-iso} implies that the exponential
  isomorphism $\exp:\ulie_\Q \to U_\Q$ is defined over $\Zp$ and hence
  over $\A$. Applying this lemma to the opposite parabolic, the
  isomorphism $\exp:\ulie_F^- \to U_{/F}^-$ is defined over $\A$ as
  well.

  Consider the subgroup schemes of $\SL_2 = \SL_{2/\Z}$:
  \begin{equation*}
    U_1^- =   \begin{pmatrix}
      1 & 0 \\
      * & 1 
    \end{pmatrix} \iso \G_{a/\Z}, \quad  T_1 = \G_{m/\Z}, \quad U_1 = 
    \begin{pmatrix}
      1 & * \\
      0 & 1 
    \end{pmatrix} \iso \G_{a/\Z}.
  \end{equation*}
  The ``big cell'' of $\SL_2$ is the subscheme $\Omega = U_1^- \cdot T_1
  \cdot U_1$; the product map defines an isomorphism $U_1^- \times T_1
  \times U_1 \to \Omega$ of schemes $/\Z$.
  
  If $\phi_{/F}:\SL_{2/F} \to G_{/F}$ is as in (1), then the restriction
  of $\phi_{/F}$ to $\Omega_{/F}$ is given by $(s,t,u) \mapsto \exp(sY_F) \cdot
  \tau_{/F}(t) \cdot \exp(uX_F)$.  Since $\tau$ and the exponential
  maps are defined over $\A$, it follows that the restriction of
  $\phi_{/F}$ to $\Omega_{/F}$ is defined over $\A$.  The proof of
  \cite[Prop. 2]{serre:PSLp} now implies that $\phi_{/F}$ is defined over
  $\A$ which settles (1). (The argument of \emph{loc.  cit.} uses that
  $\SL_{2/\A}$ is covered by $\Omega_{/\A}$ and $w\Omega_{/\A}$, for a
  suitable $w \in \SL_2(\Z)$).
  
  To prove (2), note first that
  by Lemma
  \ref{lem:rational-exponential-iso},
  $\exp:\ulie \to U$ is $P$-equivariant. Since $X_k$ is in the dense
  $P$ orbit on $\ulie$,  $\exp(X_k)$ is in the dense $P$
  orbit on $U$; since $\exp(X_k) = \phi_{/k}(
  \begin{pmatrix}
    1 & 1 \\
    0 & 1
  \end{pmatrix}),$  (a)  now follows. 
  Since by Lemma \ref{lem:rational-exponential-iso} the tangent map to
  $\exp$ is the identity, (b) holds. (c) follows by construction.
\end{proof}

\section{Why Theorem \ref{theorem:main} implies Theorems 
  \ref{theorem:testerman} and \ref{theorem:LawtherTesterman}}
\label{sec:reductions}

We return to the assumptions of the introduction; thus $k$ is
algebraically closed of characteristic $p>0$ (but not necessarily an
algebraic closure of a finite field), $G$ is reductive, and $p$ is
good for $G$.

Since Theorem \ref{theorem:LawtherTesterman} is a stronger form of
Theorem \ref{theorem:testerman}; we just prove this latter result
using Theorem \ref{theorem:main}.

\begin{proof}[Proof of Theorem \ref{theorem:LawtherTesterman}]
  First, we may replace the reductive group $G$ by its derived group,
  so that we may suppose $G$ to be semisimple.  Let $\pi:\hat G \to G$
  be the simply connected covering group of $G$. According to Lemma
  \ref{lem:central-isogeny-homeo} (in the appendix), the restrictions
  $\pi_{\mid \hat \UU}:\hat \UU \to \UU$ and $d\pi_{\mid \hat
  \NN}:\hat \NN \to \NN$ are homeomorphisms. If $u' = \pi^{-1}(u)$ and
  if $\phi:\SL_{2/k} \to \tilde G$ is a sub-principal homomorphism
  with $u'$ in its image, it follows from Remark
  \ref{rem:cochar-and-isogeny} that $\pi \circ \phi$ is a
  sub-principal homomorphism with $u$ in its image.  Thus we may
  suppose that $G$ is simply connected.
  
  It suffices to prove the theorem in the case where $k$ is an
  algebraic closure of a finite field. Indeed, let $k_0 \subset k$
  denote the algebraic closure of the prime field $\F_p$.  Then by
  \cite[Cor.  7.3]{math.RT/0007056} each nilpotent orbit in $\glie_k$
  has a point rational over $k_0$.
  
  Since $p$ is good, there is a $G$-equivariant homeomorphism $\e:\NN
  \to \UU$; see Proposition \ref{prop:springer-homeo} of the appendix.
  Moreover, by Remark \ref{rem:springer-iso-remarks} (in the
  appendix), $\e$ restricts to a $G$-equivariant homeomorphism
  $\NN_{/k_0} \to \UU_{/k_0}$.  
  This shows that also each unipotent class has a point rational over
  $k_0$.  Thus $u = gu'g^{-1}$ for some $u'$ rational over $k_0$,
  where $g$ is over $k$.  The theorem for $k_0$ and $u'$ yields a
  suitable homomorphism $\phi:SL_{2/k_0} \to G_{/k_0}$, and $\Int(g)
  \circ \phi$ then works for $k$ and $u$.
  
  Now assume that $k=k_o$, and that $u$ is a distinguished unipotent
  element.  Let $P$ be the canonical parabolic subgroup
  attached to $u$. Replacing $u$ and $P$ by a
  $G$-conjugate, we may suppose that $\Lie(P) = \plie_\Z \tensor_\Z k$
  where $\plie_\Z$ is a distinguished standard parabolic of
  $\glie_\Z$.
  
  Since $u$ has order $p$, \cite[Theorem 1]{math.RT/0007056} shows
  that a Richardson element $X$ in $\Lie(U)$ satisfies $X^{[p]} = 0$.
  Theorem \ref{theorem:main} now gives us a sub-principal homomorphism
  $\phi:\SL_{2/k} \to G_{/k}$ whose image meets the dense $P$ orbit on
  $U$.  Replacing $\phi$ by $\Int(g) \circ \phi$ for a suitable $g \in
  P$, the proof is complete for distinguished $u$.
  
  When $u$ is not distinguished, it is distinguished in a proper Levi
  subgroup $L$.  We may then apply the result in the distinguished
  case to $L$; Remark
  \ref{rem:levi-assoc-cocharacter} shows that the 
  homomorphism so obtained has the desired properties.
\end{proof}

\section{Rationality properties}

Let $F$ be a ground field, suppose that $G=G_{/F}$ is a reductive group
over $F$.

\subsection{Some rational Levi and parabolic subgroups}

We begin with a lemma:
\begin{lem}
  \label{lem:rational-parabolics}
  Let $P$ be a parabolic subgroup of $G$, and suppose that some
  maximal torus of $P$ is defined and split over $F$.  Then $P$ is
  itself defined over $F$.
\end{lem}

\begin{proof}
  Note that the our assumption means that  $G$ is
  $F$-split.  Let $T_o \le B_o$ be the ``standard'' maximal torus and
  Borel subgroup of the split group $G$. Thus $T_o$ is $F$-split, and
  representatives for the cosets in the Weyl group $W=N_G(T_o)/T_o$
  may be chosen rational over $F$.
  
  By \cite[Theorem 14.4.3]{springer98:_linear_algeb_group} one knows
  that any two $F$-split maximal tori of $G$ are conjugate by an
  element of $G(F)$. Thus we may as well suppose that $P$ contains
  $T_o$.  Choose a Borel subgroup $B$ of $P$ containing $T_o$. Since
  $B$ is also a Borel subgroup of $G$, and since all Borel subgroups
  of $G$ containing $T_o$ are conjugate by $N_G(T_o)$, we may replace
  $P$ by a $G(F)$ conjugate and suppose that $P$ contains $B_o$.  Now
  the lemma is immediate, since each of the parabolic
  subgroups of $G$ containing the standard Borel subgroup $B_o$ are
  defined over $F$.
\end{proof}

Recall \cite[\S2.9]{jantzen:Nilpotent} that when $G$ is semisimple and
the characteristic $p$ of $F$ is \emph{very good} for $G$, then the
orbit of each nilpotent and unipotent element is separable (if the
irreducible root system is different from $A_r$, then ``very good''
means the same as ``good'', while $p$ is ``very good'' for type $A_r$
provided that $r \not \congruent -1 \pmod p$).
When $G$ is reductive, we refer to \emph{loc. cit.} for
a discussion of the separability of nilpotent orbits.

\begin{theorem}
\label{theorem:nilpotent-rationality}
Let $X \in \glie(F)$ be an $F$-rational nilpotent element, and suppose
either that $F$ is perfect, or that the $G$-orbit of $X$ is separable.
  \begin{enumerate}
  \item   Then $X$ is distinguished in the Lie algebra of a Levi
    subgroup which is defined over $F$.
  \item The canonical parabolic subgroup attached to $X$ is defined
    over $F$.
  \end{enumerate}
\end{theorem}

\begin{proof}
  For the first assertion, let $C = C_G(X)$ be the centralizer of $X$.
  Under our assumptions, either $F$ is perfect, or the orbit map for
  $X$ is separable.  According to \cite[Prop.
  12.1.2]{springer98:_linear_algeb_group}, the group $C$ is then
  defined over $F$.
  
  Let $T \le C$ be a maximal torus which is defined over $F$ (such a
  torus exists by \cite[Theorem 13.3.6 and Remark
  13.3.7]{springer98:_linear_algeb_group}).  Then $L=C_G(T)$ is a Levi
  subgroup of $G$ whose Lie algebra contains $X$; see \cite[Prop.
  1.22]{DigneMichel} and \cite[Corollary
  5.4.7]{springer98:_linear_algeb_group}. The Levi subgroup $L$ is
  $F$-rational by \cite[Prop.
  13.3.1]{springer98:_linear_algeb_group}. Since $T$ is central in
  $L$, we see that the connected center of $L$ is a maximal torus of
  $C_L(X)$ hence that $X$ is distinguished in $\Lie(L)$; this proves
  (1).
  
  For (2), let $\Orbit = \Ad(G)X \subset \glie$, let $P$ be the canonical
  parabolic subgroup associated to $X$, and let $\plie =
  \Lie(P)$; see Proposition \ref{prop:canonical-parabolic}. We will
  show that $P$ (and hence also $\plie$) is defined over $F$.
  
  Fix an algebraic closure $\algc{F}$ of $F$, let $F_s$ be the
  separable closure of $F$ in $\algc{F}$, and let $\Gamma$ be the
  Galois group of $F_s/F$. Since $X$ is $F$-rational, it is stable
  under the action of $\Gamma$; since $P$ is canonically attached to
  $X$, it is also stable under the action of $\Gamma$. If $P$ is
  defined over $F_s$, (2) now follows by \cite[Proposition
  11.2.8(i)]{springer98:_linear_algeb_group}. This completes the proof
  in case $F$ is perfect. In the case where $F$ is not perfect, this
  shows that \emph{we may now suppose $F$ to be 
    separably closed}, and
  we may moreover suppose that the orbit map $G \to \Orbit$ is
  separable.
  
  Let $\mathscr{P}$ be the variety of parabolic subgroups conjugate to
  $P$, and let $P_0 \in \mathscr{P}$ denote the standard parabolic
  which is conjugate (via $G$) to $P$. Since $G$ is split over $F$,
  $P_0$ is over $F$. We identify $\mathscr{P}$ with $G/P_0$ as
  $F$-varieties. We have to show that $P$ is conjugate to $P_0$ via
  $G(F)$.
  
  Since the quotient map $G \to G/P_o$ is defined over
  $F$ and has ``local sections'' (see \cite[Lemma
  8.5.2]{springer98:_linear_algeb_group}) we may define the fiber
  space
  \begin{equation*}
    Y = G \times^{P_o} \ulie_o
  \end{equation*}
  as in \cite[Lemma 5.5.8]{springer98:_linear_algeb_group}, where
  $\ulie_o$ is the nilradical of $\Lie(P_o)$. $Y$ is an $F$-variety
  and is defined as a quotient of $G \times \ulie_o$.  Since the local
  sections of $G \to G/P_o$ are defined over $F$, a point in $Y$ is
  $F$-rational if and only if it is represented by an $F$-rational
  pair $(g,Z) \in G \times \ulie_o$. Let $p_1:Y \to \mathscr{P}$
  denote the morphism induced by $(g,Z) \mapsto \Int(g)P_o$, and let
  $p_2:Y \to \overline{\Orbit} \subset \NN$ be the morphism induced by
  $(g,Z) \mapsto \Ad(g)Z$. Then $p_1:Y \to \mathscr{P}$ is a vector
  bundle, hence $Y$ is a smooth $F$-variety.  Moreover, $p_1$ and
  $p_2$ are defined over $F$, and are equivariant (for the obvious
  action of $G$ on $Y$).  Let $\mathcal{U} = p_2^{-1}(\Orbit)$.

  Now suppose that $X$ is distinguished, so that $X$ is a Richardson
  element for $P$. Thus $\overline{\Ad(P)X} = \ulie$. Proposition
  \ref{prop:canonical-parabolic} implies that $C_G(X) = C_P(X)$.
  Since the orbit map $G \to \Orbit$ is separable, \cite[Lemma
  \S8.8]{jantzen:Nilpotent} implies that
  \begin{equation}
    \label{eq:p2-iso}
    \text{$p_2$ restricts to an isomorphism  $\mathcal{U}
    \xrightarrow{\sim} \Orbit$ of varieties.}
  \end{equation}
  
  Consider $\tilde X = p_2^{-1}(X) \in \UU$. Then $\tilde X$ is a
  simple point of $Y$ (since $Y$ is smooth).  Since $X \in \UU$,
  \eqref{eq:p2-iso} implies that $dp_2:T_{\tilde X}Y \to T_X
  \overline{\Orbit}$ is an isomorphism.  We now apply the condition
  \footnote{There is a typographic error in the statement of Corollary
    11.2.14 in \cite{springer98:_linear_algeb_group}: the condition on
    $d\phi_x$ should read ``...the tangent map $d\phi_x:T_xX \to T_yY$
    is surjective...''.}  of \cite[Corollary
  11.2.14]{springer98:_linear_algeb_group} to see that the $\tilde X$
  is rational over $F$.  Thus $\tilde X$ is represented by a pair
  $(g,Z)$ where $g \in G(F)$ and $Z \in \ulie_o(F)$, so that $P =
  gP_og^{-1}$ is defined over $F$. This proves (2) in case $X$ is
  distinguished.
  
  If $X$ is any nilpotent, then by (1) $X$ is distinguished nilpotent
  in the Lie algebra of a Levi subgroup $L$ which is defined over $F$.
  The canonical parabolic $Q$ of $X$ in the Levi subgroup $L$ is then
  defined over $F$ by the previous remarks. Recall that we suppose $F$
  to be separably closed; thus $Q$ is split over $F$.  The canonical
  parabolic subgroup $P$ of $X$ in the original group $G$ contains
  $Q$. It follows that $P$ contains a split maximal torus of $G$, so
  $P$ is defined over $F$ by Lemma \ref{lem:rational-parabolics}.
\end{proof}

\begin{theorem}
  \label{theorem:unipotent-rationality}
  Let $u \in G(F)$ be an $F$-rational unipotent element, and
  suppose either that $F$ is perfect, or that the $G$-orbit of $u$ is
  separable.
  \begin{enumerate}
  \item $u$ is distinguished in a Levi subgroup which is defined over
    $F$.
  \item The canonical parabolic
    subgroup associated to $u$ is defined over $F$.
  \end{enumerate}
\end{theorem}

\begin{proof}
  Part (1) is proved \emph{mutatis mutandum} as in part (1) of the
  previous theorem. For part (2) in the distinguished case, one must
  instead replace the vector bundle $Y$ by the ``affine-space bundle''
  $G \times^{P_o} U_o$ over $\mathscr{P}$.  The remainder of the
  argument is the same.
\end{proof}

\begin{rem}
  In the case where $F = \F_q$ is a finite field of order $q$, we can
  give a different proof of part (2) of Theorem
  \ref{theorem:nilpotent-rationality} which shows moreover that the
  \emph{grading} of \ref{sub:parabolics} is $\F_q$-rational.  Indeed,
  let $\phi:\G_m \to G$ be a cocharacter associated to $X$.  Then
  $\phi$ determines a subgroup $T_\phi = \{ (t,\phi(t)) \mid t \in
  \G_m\} \le \G_m \times G$. Since $X$ is $\F_q$-rational, the variety
  of all such subgroups $T_\phi$ (where $\phi$ ranges over all
  cocharacters associated to $X$) is stable by the action of the
  geometric Frobenius endomorphism, and is thus rational over $\F_q$
  (see Lemma \ref{lem:lie-alg-frob} and the proof of Theorem
  \ref{theorem:finite-field-case} below).  Moreover, this variety is a
  homogeneous space for the connected group $C_G^o(X)$ by Proposition
  \ref{prop:associated-cocharacter}.  Thus an application of Lang's
  Theorem \cite[Cor. 3.12]{DigneMichel} shows that there is an $\F_q$
  rational point $T_{\phi'}$.  The rationality of $T_{\phi'}$ is
  equivalent to that of $\phi'$; thus the weight spaces of $\phi'$ on
  $\glie$ are $\F_q$-rational.
\end{rem}
\begin{rem}
  Let $F$ be an arbitrary ground field, and $X$ a rational nilpotent
  element. Despite the rationality of the canonical parabolic
  associated with $X$, I do not know if there is always a cocharacter
  associated with $X$ which is defined over $F$. This is so for finite
  fields by the previous remark.
\end{rem}

\subsection{Conjugacy of nice homomorphisms}

In this section, $G$ is a reductive group over the algebraically
closed field $k$. Moreover, we suppose that $G$ is \emph{simply
connected}. Fix a distinguished parabolic subgroup $P$ of $G$, and let
$\e:\ulie \to U$ be a $P$-equivariant homeomorphism, where $U$ is the
unipotent radical of $P$ and $\ulie$ is its Lie algebra (in good
characteristic, such a homeomorphism always exists; see \S
\ref{sec:uni-vs-nil} below.)

Fix a Richardson element $u\in U$, and let $X = \e^{-1}(u)$.  We shall
say that a homomorphism $\phi:\SL_{2/k} \to G$ is \emph{nice} for $u$
with respect to $\e$ if the following property is satisfied:
\begin{equation}
\label{eq:nice}
  \phi(
  \begin{pmatrix}
    1 & 1 \\
    0 & 1 
  \end{pmatrix}) = u
 \quad \text{and} \quad
 d\phi(
  \begin{pmatrix}
    0 & 1 \\
    0 & 0
  \end{pmatrix}) = X.
\end{equation}

For a nice homomorphism $\phi$, let $\psi$ be the
co-character given by
\begin{equation}
  \label{eq:char-of-nice}
  \psi(t) = \phi(
  \begin{pmatrix}
    t & 0 \\
    0 & t^{-1}
  \end{pmatrix}), \quad t \in k^\times;
\end{equation}
since $X$ is distinguished, it is immediate that $\psi$
is associated to $X$.

\begin{prop}
  \label{prop:nice-conjugacy}
  Any two homomorphisms $\SL_{2/k} \to G$ which are nice for $u$ with
  respect to $\e$ are conjugate by an element of $C_G^o(u)$.
\end{prop}

\begin{proof}
  For $i=1,2$ let $\phi_i$ be nice homomorphisms for $u$ with respect
  to $\e$, and let $\psi_i$ be the corresponding characters as in
  \eqref{eq:char-of-nice}. According to \cite[Lemma
  5.3]{jantzen:Nilpotent}, there is an element $g \in C_G^o(u) =
  C_G^o(X)$ with $\psi_2 = \Int(g) \circ \psi_1$. Replacing $\phi_1$
  by $\Int(g) \circ \phi_1$, we may suppose that $\psi_1 = \psi_2$.
  The proposition is now a consequence of the lemma that follows.
\end{proof}

\begin{lem}
        Let $\phi_i:\SL_{2/k} \to G$,
  $i=1,2$, be nice homomorphisms for $u$ with respect to $\e$, and let
  $\psi_i$ be the corresponding cocharacters as in
  \eqref{eq:char-of-nice}.  If $\psi_1 = \psi_2$, then $\phi_1 =
  \phi_2$.
\end{lem}

\begin{proof}
  Let $\Omega \subset \SL_{2/k}$ be the big cell $\Omega =
  U_1^-T_1U_1$ as in the proof of Theorem \ref{theorem:main}. Since
  $\Omega$ is a dense subset of $\SL_{2/k}$, it suffices to show that
  the restrictions of $\phi_1$ and $\phi_2$ to $\Omega$ coincide.
  
  For $s \in k$, one has
  \begin{equation*}
    \phi_i(
    \begin{pmatrix}
      1 & s \\
      0 & 1
    \end{pmatrix}) = \Int(\psi_i(s^{1/2}))u; 
  \end{equation*}
  since $\psi_1=\psi_2$, it follows that the restrictions of the
  $\phi_i$ to $U_1$ coincide.
  
  Let $\glie(i)$ the graded components of $\glie$ with respect to
  $\psi$, and let $H = d\psi(1) \in \glie(0)$. Since $G$ is simply
  connected, Lemma \ref{lem:integral-sl2} shows that there is a unique
  $Y \in \glie(-2)$ such that $(X,H,Y)$ is an $\lie{sl}_2$ triple over
  $k$; it follows that $d\phi_i( \begin{pmatrix} 0 & 0 \\ 1 & 0
  \end{pmatrix}) = Y$ for $i=1,2$. In particular, $d\phi_1 = d\phi_2$.
  
  We may find $w \in \SL_{2/k}$  with 
  \begin{equation*}
\Ad(w)
  \begin{pmatrix}
    0 & 1 \\
    0 & 0
  \end{pmatrix} = 
  \begin{pmatrix}
    0 & 0 \\
    1 & 0
  \end{pmatrix} \quad \text{and} \quad \Int(w)
  \begin{pmatrix}
    1 & s \\
    0 & 1
  \end{pmatrix} = 
  \begin{pmatrix}
    1 & 0 \\
    s & 1
  \end{pmatrix} \quad \text{ for all $s \in k$.}     
  \end{equation*}
Since $d\phi_1 = d\phi_2$,
we find that $\phi_1(w)^{-1}\phi_2(w) \in C_G(X) = C_G(u)$.
It then follows that $\phi_1(
\begin{pmatrix}
  1 & 0 \\
  1 & 1 \\
\end{pmatrix}) = \phi_2(\begin{pmatrix}
  1 & 0 \\
  1 & 1 \\
\end{pmatrix})$. Arguing as before, one sees
that the restrictions of the $\phi_i$ to $U_1^-$ coincide, and the
lemma is proved.
\end{proof}

\begin{rem}
  \label{rem:nice-exists} Let $u \in G$ be a distinguished unipotent
  element of order $p$, and suppose that $u$ is rational over a ground
  field $F$.  We make the same hypothesis as in Theorems
  \ref{theorem:nilpotent-rationality} and
  \ref{theorem:unipotent-rationality}; thus $p$ is a good prime, and
  either the $G$-orbit of $u$ is separable, or $F$ is perfect.  Let
  $U$ be the unipotent radical of the canonical parabolic subgroup $P$
  associated with $u$ (recall by Theorem
  \ref{theorem:unipotent-rationality} that $P$ and hence $U$ are
  defined over $F$). By Lemma \ref{lem:rational-exponential-iso} there
  is a unique $P$-equivariant isomorphism $\exp:\Lie(U) \to U$ whose
  tangent map is the identity; moreover, $\exp$ is defined over $F$.
  Since $G$ is simply connected, our proof of Theorem
  \ref{theorem:LawtherTesterman} via Theorem \ref{theorem:main} shows
  that there is a sub-principal homomorphism $\phi:\SL_{2/k} \to G$
  which is nice for $u$ with respect to $\exp$. Actually, it is not
  necessary to assume simple connectivity: if $\pi:\hat G \to G$ is
  the simply connected covering group, then $\pi$ restricts to an
  isomorphism $\pi_{\mid \hat U}:\hat U \to U$ where $\hat U =
  \pi^{-1}(U)$.  By the unicity, we have $\exp \circ d\pi_{\mid
    \hat{\ulie}} = \pi_{\mid \hat{U}} \circ \hat{\exp}$, where
  $\hat{\exp}$ denotes the corresponding exponential for $\hat U$.  So
  if $\hat \phi:\SL_{2/k} \to \hat G$ is nice for $\pi^{-1}(u)$ with
  respect to $\hat{\exp}$, then $\phi = \pi \circ \hat \phi$ is nice
  for $u$ with respect to $\exp$.

  Note that there is no
  \emph{a priori} reason that $\phi$ should be defined over $F$.
\end{rem}

\subsection{Finite fields}
\label{sub:finite-fields}

Suppose now that $k$ is an algebraic closure of the finite field
$\F_q$ with $q$ elements and characteristic $p$.  Let $V$ be an affine
variety over $k$.  Recall that there is a dictionary between
$\F_q$-structures on $V$ and certain morphisms $F:V \to V$ (for
details consult e.g. \cite[Chapter 3]{DigneMichel}).  Indeed, an
$\F_q$-structure on $V$ is a finitely generated $\F_q$-subalgebra $A_0
= \F_q[V] \subset A=k[V]$ with the property that the natural map $A_o
\tensor_{\F_q} k \to A$ is an isomorphism.  The co-morphism $F^*:A \to
A$ is then given by $f \tensor \alpha \mapsto f^q \tensor \alpha$ for
$f \in A_o$ and $\alpha \in k$.  Conversely, $F$ determines $A_o$ as
$\{f \in A \mid F^* f = f^q\}$.  Note that $F^*:A \to A^q$ is
surjective, and that \cite[Prop. 3.3(i)]{DigneMichel} gives necessary and
sufficient conditions under which a surjective algebra map $A \to A^q$
determines an $\F_q$-structure on $V$. The map $F$ is called the
\emph{geometric Frobenius} endomorphism of $V$.

\begin{lem}
  \label{lem:lie-alg-frob}
  Let $H$ be a linear algebraic $k$-group defined over $\F_q$, and
 let $F$ be the corresponding Frobenius endomorphism of $H$.
 \begin{enumerate}
 \item There is a unique $q$-semilinear automorphism $\varphi$ of
   $\hlie = \Lie(H)$ such that the $\F_q$-space of
   $\varphi$-fixed-points $\hlie^\varphi$ identifies with the
   $H$-invariant $\F_q$-derivations of $\F_q[H]$.
 \item If $B \le H$ is a closed subgroup, $\Lie(F(B)) = \varphi \Lie(B)$.
 \end{enumerate}
\end{lem}

\begin{proof}
  First recall that a $\F_p$-linear automorphism $\varphi$
of a $k$-vector space  $V$ is $q$-semilinear if $\varphi(\lambda v) = \lambda^q
  \varphi(v)$ for each $v \in V$ and $\lambda \in k$.
  
  Let $A_0 = \F_q[H] \subset A = k[H]$. Recall that the
  \emph{arithmetic Frobenius} map $\varphi_a$ associated to the give
  $\F_q$-structure is the $q$-semilinear automorphism of $A$
  satisfying $\phi_a(f \tensor \alpha) = f \tensor \alpha^q$ for $f \in
  A_0$, $\alpha \in k$.
  
  Now $\hlie$ is the Lie algebra of $H$-invariant derivations
  of $A$. There is a natural map $\Der_{\F_q}(A_0) \to \Der_k(A)$. The
  $H$-invariant derivations in the image of this map form an
  $\F_q$-subspace $\hlie_0$ of $\hlie$; moreover, the natural map
  $\hlie_0 \tensor_{\F_q} k \to \hlie$ is an isomorphism. Now take for
  $\varphi$ the map satisfying $\varphi(X \tensor \alpha) = X \tensor
  \alpha^q$ (for each $X \in \hlie_0$ and $\alpha \in k$). Then
  $\hlie_0 = \hlie^\phi$ and (1) is clear.  If $X \in \hlie$ is
  regarded as a derivation, and $f \in A$, then $\varphi^{-1}(X)(f) =
  \varphi_a^{-1}(X(\varphi_a(f)))$.
  
  Let $I = \mathcal{I}(B) \normal A$ be the defining ideal of the
  closed subgroup $B$, and let $J = \mathcal{I}(FB) \normal A$ that of
  $FB$.  Since $F^*
  \circ \phi_a(f) = \phi_a \circ F^*(f) = f^q$ for each $f \in A$,
  one readily checks that $J = \phi_a(I)$.
  Thus,
  \begin{align*}
    \Lie(F(B)) &= \{X \in \hlie \mid X(f) \in J\ \forall f \in J\} \\
    &= \{X \in \hlie \mid \varphi_a^{-1}(X(\varphi_a(h))) \in I 
    \ \forall h \in I\} \\
    &= \{X \in \hlie \mid \varphi^{-1}(X)(h) = 0 \ \forall h \in I\} \\
    &= \varphi \Lie(B);
  \end{align*}
  this proves (2).
\end{proof}

We now suppose that the connected reductive group $G$ is defined over
$\F_q$ and that $p$ is good for $G$. Denote by $F$ the
corresponding Frobenius endomorphism of $G$, and by $\varphi$
the $q$-semilinear automorphism of $\glie$ as in the lemma.
Also, let $F_o$ be the Frobenius endomorphism of $\SL_{2/k}$ (and
$\varphi_0$ the $q$-semilinear automorphism of
$\lie{sl}_2(k)$) for the standard $\F_q$-structure.

\begin{theorem}
  \label{theorem:finite-field-case}
  Let $u \in G$ be an $\F_q$-rational unipotent element of order $p$.
  Then there is a sub-principal homomorphism $\psi:\SL_{2/\F_q} \to
  G_{/\F_q}$ defined over $\F_q$ whose image contains $u$.
\end{theorem}

\begin{proof}
  As in the proof of Theorem \ref{theorem:LawtherTesterman}, we may
  suppose that $G$ is simply connected (note that the simply connected
  covering $\pi:\hat G \to G$ is defined over $\F_q$. In fact, for
  \emph{any} ground field $F$, the simply connected cover of an
  $F$-reductive group is again over $F$; see
  \cite[Lemma 16.2.4]{springer98:_linear_algeb_group}.)

  By Theorem \ref{theorem:unipotent-rationality}(1), we may suppose
  that $u$ is distinguished.  Let $P$ be the canonical parabolic
  associated to $u$, and let $U$ be its unipotent radical (by Theorem
  \ref{theorem:unipotent-rationality} these subgroups are defined over
  $\F_q$). By Remark \ref{rem:nice-exists}, there is a sub-principal
  homomorphism $\psi$ which is nice for $u$ with respect to $\exp$,
  where $\exp:\Lie(U) \to U$ is the exponential $\F_q$-isomorphism of
  Lemma \ref{lem:rational-exponential-iso}.  For any sub-principal
  homomorphism $\psi$ nice for $u$ with respect to $\exp$, we get a
  subgroup \begin{equation*} \Gamma = \Gamma_\psi = \{(g,\psi(g)) \mid
  g \in \SL_{2/k} \} \le \SL_{2/k} \times G \end{equation*}
  satisfying:
\begin{enumerate} \item $\Gamma$ contains $u'=( \begin{pmatrix} 1 & 1
\\ 0 & 1 \end{pmatrix},u)$.  \item $\Lie(\Gamma)$ contains $X'=(
\begin{pmatrix} 0 & 1 \\ 0 & 0 \end{pmatrix}, X)$ (where $u =
\exp(X)$).  \item The restriction of the first projection $p_1$ to
$\Gamma$ is an isomorphism (of algebraic groups).  \end{enumerate} On
the other hand, if $\Gamma \le \SL_{2/k} \times G$ is any subgroup
satisfying (1), (2), and (3), then $\Gamma$ is the graph of a unique
sub-principal homomorphism $\psi:\SL_{2/k} \to G$ which is nice for
$u$ with respect to $\exp$.
  
  Since $G$ is simply connected, it follows from Proposition
  \ref{prop:nice-conjugacy} that the variety $\mathcal{V}$ of all
  subgroups $\Gamma$ satisfying (1), (2), and (3) is a homogeneous
  space for the connected group $C_G^o(u)$.
  
  We now claim that the variety $\mathcal{V}$ is stable by the
  Frobenius $F_1=F_o \times F$ of $\SL_{2/k} \times G$.  Let $\Gamma$
  in $\mathcal{V}$; we verify that (1)--(3) hold for $F_1\Gamma$.
  Since $u'$ is clearly $F_1$ stable, (1) holds.  Since $\exp$ is
  defined over $\F_q$, $X' \in \lie{sl}_2(k) \oplus \glie$ is fixed by
  $\varphi_1=\varphi_0 \oplus \varphi$.  According to the lemma, the
  Lie algebra of $F_1 \Gamma$ is $\varphi_1\Lie(\Gamma)$; it follows
  that $X' \in \Lie(F_1 \Gamma)$ so that (2) holds.
  
  To verify (3), note that the restriction of $p_1$ to $F_1 \Gamma$ is
  evidently bijective (since $F_1$ itself is bijective). It remains to
  see that the restriction of the differential $dp_1$ to $\Lie(F_1
  \Gamma)$ is bijective; that follows immediately from the equality
  $\Lie(F_1 \Gamma) = \varphi\Lie(\Gamma)$ proved in the lemma.
  
  We have now verified that $\mathcal{V}$ is $F_1$ stable. Since this
  variety is a homogeneous space for the connected group $C_G^o(u)$,
  an application of Lang's Theorem \cite[Cor.  3.12]{DigneMichel}
  yields a point $\Gamma \in \mathcal{V}$ fixed by $F_1$.  The
  homomorphism $\psi$ whose graph is $\Gamma$ then has the desired
  properties.
\end{proof}

\begin{rem}
  The theorem yields in particular a homomorphism
  $\SL_2(\F_q) \to G(\F_q)$ between the groups of rational points.
\end{rem}

\begin{rem}  
  It is common in the study of finite simple groups to consider a more
  general notion of Frobenius endomorphism. Let $G$ be a semisimple
  group over $k$. A surjective endomorphism $\sigma$ of $G$ will be
  called a Frobenius endomorphism provided the fixed-point group
  $G^\sigma$ is finite; such endomorphisms are thoroughly studied by
  Steinberg in \cite{steinberg-endomorphisms}.
  
  Let $\sigma$ be a Frobenius endomorphism of the semisimple group
  $G$, and suppose that $u \in G$ is $\sigma$-stable and of order $p$.
  It is proved in \cite[Theorem 5.1]{Proud-Saxl-Testerman} that in
  case $\sigma$ is a $q$-Frobenius endomorphism (for the definition,
  see \cite[\S 1]{Proud-Saxl-Testerman}), there is a $\sigma$-stable
  closed $A_1$-type subgroup $S \le G$ containing $u$.  In general,
  there need be no such subgroup; see \cite[Lemma
  2.1(i) and Lemma 2.2]{Proud-Saxl-Testerman}.
  
  We observe that $\sigma$ is a $q$-Frobenius endomorphism (in the
  sense of \cite{Proud-Saxl-Testerman}) if and only if it is the
  geometric Frobenius endomorphism associated to an $\F_q$-structure
  of $G$.  
  This observation was made in a slightly different context in 11.6
  (p. 76) of Steinberg's paper on Endomorphisms of Linear Algebraic
  Groups \cite{steinberg-endomorphisms}; note that if $\sigma$ is a
  $q$-Frobenius, the results in \cite[Chapter 3]{DigneMichel} (cited
  above) yield the existence of a suitable $\F_q$-structure, while if
  $\sigma$ is not a $q$-Frobenius, then $\sigma^*A$ is not equal to
  $A^{q'}$ for any $q' = p^a$, where $A = k[G]$.
  
  Thus, when $G$ is semisimple in good characteristic, the content of
  Theorem \ref{theorem:finite-field-case} is roughly that of
  \cite[Theorem 5.1]{Proud-Saxl-Testerman}. Note however that on the
  one hand, the theorem in \cite[Theorem 5.1]{Proud-Saxl-Testerman}
  treats also bad primes (the reader is referred there for
  a precise statement in bad characteristic), while on the other hand,
  Theorem \ref{theorem:finite-field-case} gives more precise
  information about the $A_1$-embedding.
  
  Since the proof of Theorem \ref{theorem:finite-field-case} is
  achieved essentially through an application of Lang's theorem, and
  since Lang's theorem remains valid for any Frobenius endomorphism
  $\sigma$, the reader may be curious why Theorem
  \ref{theorem:finite-field-case} is not valid for an arbitrary
  $\sigma$. The reason is as follows: If $\mathcal{V}$ is the variety
  appearing in the proof, and if $\sigma$ is not the
  geometric Frobenius for an $\F_q$-structure on $G$, then
  $\mathcal{V}$ need not be $(F_0 \times \sigma)$-stable (note that
  there is no analogue of Lemma \ref{lem:lie-alg-frob} for $\sigma$).
\end{rem}

\section{Appendix: Comparing the unipotent and nilpotent varieties}
\label{sec:uni-vs-nil}

Let $k$ be an algebraically closed field with characteristic $p \ge 0$,
and $G$ a reductive group over $k$.

\begin{lem}
  \label{lem:p-power-torus}
  Suppose $p>0$, and consider an algebraic torus $T$ over $k$.  Let
  $\zlie \subset \Lie(T)$ be a $p$-Lie subalgebra.  Then $A \mapsto
  A^{[p^n]}$ defines a homeomorphism $\zlie \to \zlie$ for each $n \ge
  1$.
\end{lem}

\begin{proof}
  We first observe that the map $\varphi:\Aff^m_{/k} \to \Aff^m_{/k}$
  given by
  \begin{equation*}
    (x_1,\dots,x_m) \mapsto (x_1^{p^n},\dots,x_m^{p^n})
  \end{equation*}
  is a homeomorphism for all $m,n \ge 1$. Indeed, this map is a
  morphism of varieties, hence continuous; moreover, it is evidently
  bijective. That it is a homeomorphism will follow provided that it
  is open.  If $g$ is a regular (i.e. polynomial) function on
  $\Aff^m_{/k}$, let $D(g)$ denote the distinguished open subset of
  $\Aff^m_{/k}$ defined by it.  There is a polynomial function $f$
  with $g^{p^n} = \varphi^*f$; thus $\varphi(D(g)) =
  \varphi(D(g^{p^n})) = D(f)$ is open as desired.
  
  Let $\sigma:\Lie(T) \to \Lie(T)$ be the map $A \mapsto A^{[p]}$; for
  a subspace $\zlie \subset \Lie(T)$, we write $\zlie^\sigma = \{A \in
  \zlie \mid A = \sigma(A)\}$. Since $\Lie(T)$ is an Abelian algebra,
  $\sigma$ is additive and ``semilinear'': $\sigma(\alpha A) =
  \alpha^p \sigma(A)$ for $\alpha \in k$ and $A \in \Lie(T)$. Thus
  $\zlie^\sigma$ is an $\F_p$-vector space.
  
  One knows that the canonical map $\Lie(T)^\sigma \tensor_{\F_p} k
  \to \Lie(T)$ is an isomorphism (indeed: it suffices to observe that
  this is true when $T = \G_m$).
  
  We now claim that a $k$-subspace $\zlie \subset \Lie(T)$ is a
  $p$-subalgebra if and only if the canonical map $\zlie^\sigma
  \tensor_{\F_p} k \to \zlie$ is an isomorphism.
  
  This claim follows from the (apparently) more general statement:
  suppose that $V$ is a finite dimensional $k$-vector space, that
  $\sigma:V \to V$ is a bijective, additive, semilinear map, and that
  $V^\sigma \tensor_{\F_p} k \to V$ is an isomorphism. Then a
  $k$-subspace $W \subset V$ is $\sigma$ stable if and only if the
  canonical map $W^\sigma \tensor_{\F_p} k \to W$ is an isomorphism;
  for this, see the proof of \cite[Proposition AG 14.2]{Bor1}.
  
  To finish the proof of the lemma, note that the choice of an
  $\F_p$-basis for $\zlie^\sigma$ identifies the map $(A \mapsto
  \sigma^n(A) = A^{[p^n]}):\zlie \to \zlie$ with $\varphi:\Aff^{\dim
    \zlie}_{/k} \to \Aff^{\dim \zlie}_{/k}$ and the lemma follows.
\end{proof}

In the situation of the lemma, we will write $\Theta_n:\zlie \to
\zlie$ for the inverse of the homeomorphism $A \mapsto A^{[p^n]}$, for
$n \ge 1$.  Thus $\Theta_n(A)$ is a sort of ``$p^n$-th root'' of $A
\in \zlie$.  Note that $\Theta_n$ is not a morphism of varieties
(since the morphism $A \mapsto A^{[p^n]}$ is purely inseparable of degree
$p^n$).

We denote by $\UU(G) =
\UU$ the variety of unipotent elements in $G$, and by
$\NN(G) = \NN$ the variety of nilpotent elements in $\glie$.

\begin{lem}
  \label{lem:central-isogeny-homeo}
  Let $\pi:\hat G \to G$ be a central isogeny of reductive groups over
  $k$. Then $\pi$ restricts to a homeomorphism $\pi_{\mid \hat
    \UU}:\hat \UU \xrightarrow{\sim} \UU$, and $d\pi$ restricts to a
  homeomorphism $d\pi_{\mid \hat \NN}:\hat \NN \xrightarrow{\sim}
  \NN$.  If the characteristic of $k$ is 0 or if $d\pi:\hat \glie \to
  \glie$ is bijective, these maps are isomorphisms of varieties.
\end{lem}

\begin{proof}
  Let us recall from \cite[\S 22]{Bor1} that ``$\pi$ is a central isogeny''
  means that $Z=\ker \pi$ is finite and hence central (we will regard
  this kernel as a (reduced) group variety rather than as a group
  scheme), and that $\zlie = \ker d\pi$ is central.
  
  First, we note that $\pi_{\mid \hat \UU}$ and $d\pi_{\mid \hat \NN}$
  are bijective; for the latter, this is proved in \cite[Prop.
  2.6]{jantzen:Nilpotent}; the argument for $\pi_{\mid \hat \UU}$ is
  the same.
  
  If $p=0$, then $d\pi:\hat \glie \to \glie$ is an
  isomorphism (since the kernel of $\pi$ is finite).  If $d\pi:\hat
  \glie \to \glie$ is an isomorphism, then $\pi_{\mid \hat \UU}$ and
  $d\pi_{\mid \hat \NN}$ are isomorphisms of varieties by
  \cite[Theorem 5.3.2]{springer98:_linear_algeb_group}.  Thus, we may
  now suppose that $p>0$.

  It remains to show that $\pi_{\mid \hat\UU}$ and $d\pi_{\mid \hat
    \NN}$ are open maps.  Note first that $(*)$ if $f:X \to Y$ is an
  open map between topological spaces, and $X' \subset
  X$ satisfies $X' = f^{-1}(f(X'))$ then $f_{\mid X'}:X' \to f(X')$ is
  open.
  
  Any surjective morphism of algebraic groups is open; this follows
  from \cite[Theorem 5.1.6(i)]{springer98:_linear_algeb_group}. In
  particular, $\pi:\hat G \to G$ is open, and $d\pi:\hat \glie \to
  d\pi(\hat \glie)$ is open.
  
  Let $\VV = \pi^{-1}(\UU)$. Then $(*)$ shows that $\pi_{\mid \VV}:\VV
  \to \UU$ is an open map. Since $Z$ is a finite subgroup of $G$, the
  connected components of the variety $\VV$ are precisely the sets
  $z\hat \UU$ where $z \in Z$. In particular, $\hat \UU$ is an open
  subset of $\VV$. This shows that $\pi_{\mid \hat \UU}$ is open as
  desired.
  
  Finally, let $\MM = d\pi^{-1}(\NN)$. Since $\NN \subset d\pi(\hat
  \glie)$, $(*)$ shows that $d\pi_{\mid \MM}:\MM \to \NN$ is an open map.
  We claim that the map
  \begin{equation*}
    \Phi:\zlie \times \hat \NN \to \MM \quad \text{via} \ (Z,X) \mapsto Z + X
  \end{equation*}
  is a homeomorphism. This map is a morphism hence continuous.  We
  will produce an explicit inverse. Let $n \ge 1$ have the property
  that $X^{[p^n]} = 0$ for all $X \in \hat \NN$ [it suffices to choose
  $n$ so that $X^{[p^n]} = 0$ for a \emph{regular} nilpotent element
  $X$]. For any $(Z,X) \in \zlie \times \hat \NN$, we have $[Z,X] = 0$
  so that $(Z+X)^{[p^n]} = Z^{[p^n]} \in \zlie$.  Denoting by
  $\Theta_n:\zlie \to \zlie$ the ``$p^n$-th root'' map as in the
  remarks following the previous lemma, we may define
  \begin{equation*}
    \Psi:\MM \to \zlie \times \hat \NN \quad \text{via} \ A 
    \mapsto (\Theta_n(A^{[p^n]}),A - \Theta_n(A^{[p^n]})).
  \end{equation*}
  Then $\Psi$ is continuous by construction, and $\Phi$ and $\Psi$ are
  inverse homeomorphisms.
  
  It now follows that the map $\Gamma:\zlie \times \hat \NN \to \NN$
  given by $(Z,X) \mapsto d\pi(X)$ is open (since $\Gamma = d\pi_{\mid
    \MM} \circ \Phi$). If $\mathcal{W} \subset \hat \NN$ is an open
  set, then $\zlie \times \mathcal{W} \subset \zlie \times \hat \NN$
  is an open set.  Thus $\Gamma(\zlie \times \mathcal{W}) =
  d\pi(\mathcal{W})$ is open. We have showed that $d\pi_{\mid \hat
    \NN}$ is an open map, which completes the proof of the lemma.
\end{proof}

\begin{lem}
  \label{lem:iso-respects-parabolics}
  Let $\e:\NN \to \UU$ be a $G$-equivariant homeomorphism.  Let $P
  \subset G$ be a parabolic subgroup with Levi decomposition $P =
  L\cdot V$.  Then $\e$ restricts to a $P$-equivariant
  homeomorphism $\vlie = \Lie(V) \to V$, and to an $L$-equivariant
  homeomorphism $\NN(L) \to \UU(L)$.
\end{lem}

\begin{proof}
  Suppose first that $P=B$ is a Borel subgroup with unipotent radical
  $U$. Then the proof given in \cite[Proof of Theorem 5.9.6 (2nd
  paragraph)]{Carter1} shows that $\e$ induces a homeomorphism $\ulie
  = \Lie(U) \to U$, where $U$ is the unipotent radical of $B$ (in
  \emph{loc. cit.} one is in the situation where $\e$ is an
  isomorphism of varieties, but the argument depends only on
  topological properties of $\e$).
  
  Now suppose that $P$ is a parabolic subgroup containing $B$, and
  that $P = L\cdot U_P$ is a Levi decomposition.  Let $\NN(L)$ and
  $\UU(L)$ denote the nilpotent and unipotent varieties of $L$
  (regarded as subvarieties of $\NN$ and $\UU$). Let $U^-$ denote the
  unipotent radical of the Borel group opposite to $B$. Then we have:
  \begin{gather*}
    \UU(L) \cap U = \{u \in U \mid \Int(L)u \cap U^- \not = \emptyset\}, \\
    \NN(L) \cap \ulie = \{X \in \ulie \mid \Ad(L)X \cap 
    \ulie^- \not = \emptyset \}
  \end{gather*}
  and
  \begin{equation*}
    U_P = \{u \in U \mid \Int(P)u \subset U\},
    \quad \Lie(U_P) = \{X \in \ulie \mid \Ad(P)X \subset \ulie\}.
  \end{equation*}
  The required properties of $\e$ are now immediate from equivariance.
\end{proof}

\begin{prop}
  \label{prop:springer-homeo}
  Suppose that $p$ is good for $G$. Then there
  is a $G$-equivariant homeomorphism $\e:\NN \to \UU$.
\end{prop}

\begin{proof}
  Suppose first that $G$ is simply connected and semisimple. Then the
  result is due to Springer; see \cite[Theorem 6.20]{hum-conjugacy}.
  In fact, one gets in this case an isomorphism of varieties; see
  \cite[Cor. 9.3.4]{Bard-Rich-LunaSlice}.
  
  One now deduces the result when $G$ is the product of a torus and a
  simply connected semisimple group.  Since there is a central isogeny
  from such a group onto our reductive group $G$ \cite[Theorem
  9.6.5]{springer98:_linear_algeb_group}, the result follows from
  Lemma \ref{lem:central-isogeny-homeo}.
\end{proof}

\begin{rem}
 \label{rem:springer-iso-remarks}
 Let $k_0 \subset k$ be a field extension with both $k_0$ and $k$
 algebraically closed. Then the homeomorphism $\e:\NN_{/k} \to
 \UU_{/k}$ of the proposition may be chosen so that its restriction to
 $k_0$ points defines a homeomorphism $\NN_{/k_0} \to \UU_{/k_0}$.
 Indeed, in the case where $G$ is simply connected, there is an
 equivariant isomorphism between the two varieties which is defined
 over $\Z[1/f]$, where $f$ is the product of the bad primes, and hence
 over $k_0$; see \cite[\S 6.21]{hum-conjugacy}.  Thus the claim is
 true in the simply connected case. For the general statement, there
 is a central $k_0$-isogeny from a simply connected group to $G$. 
The homeomorphisms in Lemma
 \ref{lem:central-isogeny-homeo} are then defined over $k_0$, whence
the result in general.
%, the homeomorphisms $d\pi_{\mid \hat
%   \NN_{/k}}:\hat \NN_{/k} \to \NN_{/k}$ and $\pi_{\mid \hat
%   \UU_{/k}}:\hat \UU_{/k} \to \UU_{/k}$ are morphisms of varieties
% defined over $k_0$ which restrict to homeomorphisms on $k_0$ points.;
\end{rem}

\section{Appendix: Springer's isomorphism and $p$-th powers}

Let $k$ and $G$ be as in the previous appendix, and assume that $p>0$.
Denote by $\UU$ the variety of unipotent elements in $G$, and let
$\NN$ be the variety of nilpotent elements in $\glie$.  We suppose the
following hypothesis to hold:
\begin{itemize}
\item[$(*)$] the group $G$ has a faithful rational representation
  $(\rho,V)$ for which the trace form $\beta(X,Y) = tr(d\rho(X)\circ
  d\rho(Y))$ on $\glie$ is nondegenerate.
\end{itemize}

For convenience, we identify $G$ with a subgroup of $\GL(V)$, and
hence also $\glie$ with a subalgebra of $\lie{gl}(V)$. Thus the
$p$-power map $X \mapsto X^{[p]}$ on the $p$-Lie algebra $\glie$ is
the restriction of the usual $p$-power map in the associative algebra
$\End(V)$. The hypothesis $(*)$ implies that $\lie{gl}(V) = \mm \oplus
\glie$, where $\mm = \glie^\perp = \{X \in \lie{gl}(V) \mid \beta(X,Y)
= 0 \ \text{for each}\ Y \in \glie\}$.

Let $T \subset G$ be a maximal torus, and let $T_1$ be any maximal
torus of $\GL(V)$ containing $T$. Write $\hlie$ for the Lie algebra of
$T$, and $\hlie_1$ for that of $T_1$.
\begin{lem}
\label{lem:torus-orthog}
  $\hlie_1 = (\hlie_1 \cap \mm) \oplus
  \hlie$. In particular, $\hlie_1 \cap \mm$ is the orthogonal
  complement of $\hlie$ in $\hlie_1$ with respect to $\beta$.
\end{lem}

\begin{proof}
  Since the restriction of $\beta$ to $\hlie$ is non-degenerate, the
  second assertion is a consequence of the first.  Let $Y \in
  \hlie_1$, and write $Y = A + B$ with $A \in \glie$ and $B \in \mm$.
  The first assertion follows if we show that $A \in \hlie$.  For each
  weight $\lambda \in X^*(T)$, we have a $T$-module decomposition
  $\lie{gl}(V)_\lambda = \glie_\lambda \oplus \mm_\lambda$. This
  applies in particular for $\lambda = 0$; since $\hlie = \glie_0$ we
  get $\lie{gl}(V)_0 = \hlie \oplus \mm_0$.  As $T$ acts trivially on
  $\hlie_1$, we have $Y \in \lie{gl}(V)_0$ so that $A \in \hlie$ as
  desired.
\end{proof}

\begin{lem}
\label{lem:torus-p-power}
  For each $Y \in \hlie$, there is an element $Z \in \hlie$ 
  with $Z^{[p]} = Y$.
\end{lem}

\begin{proof}
  It suffices to show that $\hlie$ has a $k$-basis $\{H_i\}$
  consisting of elements with $H_i^{[p]} = H_i$; indeed, if there is
  such a basis, and if $Y = \sum_i a_i H_i$, then $Z = \sum_i
  a_i^{1/p} H_i$ works (since $\hlie$ is Abelian). To see that $\hlie$
  has such a basis, choose an isomorphism $T \iso (\G_m)^n$ and use
  the fact that $\Lie(\G_m)$ has for basis element the
  $\G_m$-invariant derivation $H = T^{-1}\frac{d}{dT}$ of $k[\G_m] =
  k[T^{\pm 1}]$, which satisfies $H^{[p]} = H$.
\end{proof}

\begin{lem}
\label{lem:perp-p-power}
The subspace $\mm \cap \hlie'$ is invariant by the $p$-power map: i.e.
if $Y \in \mm \cap \hlie'$, then $Y^{[p]} \in \mm \cap \hlie'$.
\end{lem}

\begin{proof}
  We know that $\mm \cap \hlie' = \{Y \in \hlie' \mid \beta(Y,H) = 0
  \text{ for each } H \in \hlie\}$. Let $Y \in \mm \cap \hlie'$.  The
  lemma then follows if we show that $\beta(Y^{[p]},H) = 0$ for each
  $H \in \hlie$. Using the previous lemma, we may find an element $Z
  \in \hlie$ with $Z^{[p]} = H$. Then we have $\beta(Y^{[p]},H) =
  \beta(Y^{[p]},Z^{[p]}) = \beta(Y,Z)^p = 0$, where the second
  equality holds since $\beta$ is the \emph{trace form} of the
  representation $(\rho,V)$. This shows that $\beta(Y^{[p]},H) = 0$, as
  desired.
\end{proof}

Recall that we regard $G$ as a subgroup of $\GL(V)$, and hence as a
subset of $\End(V) = \lie{gl}(V)$.
\begin{lem}
  Let $\pi:G \to \glie$ be the restriction of the orthogonal
  projection with respect to $\beta$. Then $\pi(g^p) = \pi(g)^{[p]}$
  for each $g \in G$.
\end{lem}

\begin{proof}
  Since $g \mapsto \pi(g^p)$ and $g \mapsto \pi(g)^{[p]}$ are
  morphisms of algebraic varieties, it is enough to show 
  that they coincide on a dense subset of $G$. The semisimple elements
  in $G$ constitute such a dense set, by \cite[Theorem
  2.5]{hum-conjugacy}. If $s \in G$ is semisimple, then $s$ lies in a
  maximal torus $T$ of $G$ by \cite[Theorem
  6.3.5(i)]{springer98:_linear_algeb_group}. Choose a maximal torus
  $T_1 \le \GL(V)$ containing $T$, and write $\hlie$, $\hlie_1$ for
  their Lie algebras as before. Identifying the Lie algebra of
  $\GL(V)$ with $\End(V)$, we may regard the torus $T$ as a subset of
  $\hlie_1$. In particular, we may regard $s$ as an element of
  $\hlie_1$. According to Lemma \ref{lem:torus-orthog} we may write $s
  = A + H$ with $A \in \mm \cap \hlie_1$ and $H \in \hlie$.  Since
  $\hlie_1$ is an Abelian Lie algebra, we have $s^p = s^{[p]} =
  A^{[p]} + H^{[p]}$. According to  Lemma \ref{lem:perp-p-power} the
  subspace $\mm \cap \hlie_1$ is closed under $p$-powers; thus
  $A^{[p]} \in \mm$. It follows that $\pi(s^p) = H^{[p]} =
  \pi(s)^{[p]}$ as desired.
\end{proof}

\begin{theorem}
  Suppose that $G$ is quasisimple, that $p$ is good for $G$, and
  that either $G = GL_n$ or $G$ is almost simple and its root system
  is not of type $A_n$. Then there is a $G$-equivariant isomorphism of
  varieties $L:\UU \to \NN$ such that $L(x^p) = L(x)^{[p]}$ for all $x
  \in \UU$.
\end{theorem}

\begin{proof}
  The existence of a $G$-morphism $L$ (without the condition on $p$-th
  powers) follows from \cite{Bard-Rich-LunaSlice}.  The construction
  in \emph{loc. cit.}  proceeds as follows.  If $G = \GL_n$, one takes
  for $L$ the map $x \mapsto x-1$ (and the condition on $p$-powers is
  clear in that case). Otherwise, according to \cite[Ch. 1 Lemma
  5.3]{springer-steinberg}, a group isogenous to $G$ satisfies the
  hypothesis $(*)$ and such that the identity endomorphism of the
  representation $V$ is orthogonal via the trace form $\beta$ to every
  element of $\glie$ (equivalently: each element of $\glie$ has trace
  0). According to the summary in \cite[\S 0.13]{hum-conjugacy},
  $\glie$ is a simple Lie algebra under our assumptions. If $\tilde G$
  is a group isogenous to $G$, the isogeny thus induces an isomorphism
  on Lie algebras. Since the isomorphisms of Lemma
  \ref{lem:central-isogeny-homeo} evidently respect the $p$-power
  operations, we may if necessary replace $G$ with an isogenous group
  and so suppose that $G$ itself satisfies
  (*).
  
  One then considers the map $\pi:G \to \glie$ which is the
  restriction of the orthogonal projection with respect to the trace
  form. Then $\pi$ is clearly $G$-equivariant, and satisfies $\pi(1) =
  0$. According to \cite{Bard-Rich-LunaSlice}, $L = \pi_{\mid \UU}:\UU
  \to \NN$ is a $G$-isomorphism of varieties.  It follows from the
  previous lemma that $L(x^p) = L(x)^{[p]}$, whence the theorem.
\end{proof}

\begin{rem}
  This theorem permits a simplification of the proof of Testerman's
  ``order formula'' given in \cite{math.RT/0007056}. It allows one to
  deduce the order formula for unipotent elements from the
  corresponding formula for the $p$-nilpotence degree of nilpotent
  elements; the comparison of these respective ``orders'' was achieved
  in a different way in \emph{loc. cit.}
\end{rem}

%\bibliography{MathBib}

\begin{thebibliography}{Jantzen}

\bibitem[Bor70]{borel70:_proper_linear_repres_cheval_group}
Armand Borel, \emph{Properties and linear representations of {C}hevalley
  groups}, Seminar on Algebraic Groups and Related Finite Groups (The Institute
  for Advanced Study, Princeton, N.J., 1968/69), Springer, Berlin, 1970,
  Lecture Notes in Mathematics, Vol. 131, pp.~1--55.

\bibitem[Bor91]{Bor1}
\bysame, \emph{Linear algebraic groups}, 2nd ed., Grad. Texts in Math.,
  no. 129, Springer Verlag, 1991.

%\bibitem[Bou89]{BouLie123}
%N.~Bourbaki, \emph{{L}ie groups and {L}ie algebras, chapters 1,2,3},
%  Springer-Verlag, Berlin, 1989.

\bibitem[BR85]{Bard-Rich-LunaSlice}
Peter Bardsley and R.~W. Richardson, \emph{\'{E}tale slices for algebraic
  transformation groups in characteristic $p$}, Proc. London Math. Soc. (3)
  \textbf{51} (1985), no.~2, 295--317.

\bibitem[Car93]{Carter1}
Roger~W. Carter, \emph{Finite groups of {L}ie type: conjugacy classes and
  complex characters}, John Wiley \& Sons Ltd., Chichester, 1993, Reprint of
  the 1985 original.

\bibitem[DM91]{DigneMichel}
F.~Digne and J.~Michel, \emph{Representations of finite groups of {L}ie type},
  London Math. Soc. Student Texts, vol.~21, Cambridge University Press,
  Cambridge, 1991.

\bibitem[Hu95]{hum-conjugacy}
James~E. Humphreys, \emph{Conjugacy classes in semisimple algebraic groups},
  Math. Surveys and Monographs, vol.~43, Amer. Math. Soc., 1995.

\bibitem[J]{jantzen:Nilpotent}
Jens~Carsten Jantzen, \emph{Nilpotent orbits in representation theory}, 
  Notes from Odense summer school, August 2000.

\bibitem[LT99]{Lawther-Testerman}
R.~Lawther and D.~M. Testerman, \emph{${A}\sb 1$ subgroups of exceptional
  algebraic groups}, Mem. Amer. Math. Soc. \textbf{141} (1999), no.~674,
  viii+131.

\bibitem[M02]{math.RT/0007056}
George~J. McNinch, \emph{Abelian unipotent subgroups of reductive groups}, J.
  Pure Appl. Algebra \textbf{167} (2002), 269--300,
  \mbox{arXiv:math.RT/0007056}.


\bibitem[Pom77]{PommereningI}
Klaus Pommerening, \emph{\"{U}ber die unipotenten {K}lassen reduktiver
  {G}ruppen}, J. Algebra \textbf{49} (1977), no.~2, 525--536.

\bibitem[Pom80]{PommereningII}
\bysame, \emph{\"{U}ber die unipotenten {K}lassen reduktiver
  {G}ruppen. {I}{I}}, J. Algebra \textbf{65} (1980), no.~2, 373--398.

\bibitem[PST00]{Proud-Saxl-Testerman}
Richard Proud, Jan Saxl, and Donna Testerman, \emph{Subgroups of type ${A}\sb
  1$ containing a fixed unipotent element in an algebraic group}, J. Algebra
  \textbf{231} (2000), no.~1, 53--66. 

\bibitem[Sei00]{seitz-unipotent}
Gary~M. Seitz, \emph{Unipotent elements, tilting modules, and saturation},
  Invent. Math \textbf{141} (2000), 467--502.

\bibitem[Ser96]{serre:PSLp}
Jean-Pierre Serre, \emph{Exemples de plongements des groupes {${\text{PSL}}\sb
  2({\mathbf {F}}\sb p)$} dans des groupes de {L}ie simples}, Invent. Math.
  \textbf{124} (1996), no.~1-3, 525--562.

\bibitem[Spa84]{spalt-transverse}
Nicolas Spaltenstein, \emph{Existence of good transversal slices to nilpotent
  orbits in good characteristic}, J. Fac. Sci. Univ. Tokyo Sect. IA Math.
  \textbf{31} (1984), no.~2, 283--286.

\bibitem[Spr98]{springer98:_linear_algeb_group}
Tonny~A. Springer, \emph{Linear algebraic groups}, 2nd ed., Progr. in Math.,
  vol.~9, Birkh{\"a}user, Boston, 1998.

\bibitem[SS70]{springer-steinberg}
Tonny~A. Springer and Robert Steinberg, \emph{Conjugacy classes}, Seminar on
  Algebraic Groups and Related Finite Groups (The Institute for Advanced Study,
  Princeton, N.J., 1968/69), Springer, Berlin, 1970, Lecture Notes in
  Mathematics, Vol. 131, pp.~167--266.

\bibitem[Ste68]{Steinberg}
Robert Steinberg, \emph{Lectures on {C}hevalley groups}, Yale University, 1968.

\bibitem[Ste68a]{steinberg-endomorphisms}
\bysame, \emph{Endomorphisms of linear algebraic groups}, American
  Mathematical Society, Providence, R.I., 1968.

\bibitem[Tes95]{testerman}
Donna Testerman, \emph{{$A_1$}-type overgroups of elements of order $p$ in
  semisimple algebraic groups and the associated finite groups}, J. Algebra
  \textbf{177} (1995), 34--76.

\end{thebibliography}
\providecommand{\bysame}{\leavevmode\hbox to3em{\hrulefill}\thinspace}

\end{document}